\documentclass[11pt]{amsart}%

\usepackage[all]{xy}
\usepackage{amsfonts,amssymb,amsmath,amscd}
\usepackage{amsthm}

\newcounter{zlist}
\newenvironment{zlist}{\begin{list}{{\rm(\arabic{zlist})}}{
\usecounter{zlist}\leftmargin2.5em\labelwidth2em\labelsep0.5em
\topsep0.6ex\itemsep0.3ex plus0.2ex minus0.3ex
\parsep0.3ex plus0.2ex minus0.1ex}}{\end{list}}

\newcounter{blist}
\newenvironment{blist}{\begin{list}{{\rm(\alph{blist})}}{
\usecounter{blist}\leftmargin2.5em\labelwidth2em\labelsep0.5em
\topsep0.6ex \itemsep0.3ex plus0.2ex minus0.3ex
\parsep0.3ex plus0.2ex minus0.1ex}}{\end{list}}

\newcounter{rlist}
\newenvironment{rlist}{\begin{list}{{\rm(\roman{rlist})}}{
\usecounter{rlist}\leftmargin2.5em\labelwidth2em\labelsep0.5em
\topsep0.6ex\itemsep0.3ex plus0.2ex minus0.3ex
\parsep0.3ex plus0.2ex minus0.1ex}}{\end{list}}

\newenvironment{bull}{\begin{list}{{$\bullet$ }}{
 \leftmargin2.5em\labelwidth2em\labelsep0.5em
\topsep0.6ex\itemsep0.3ex plus0.2ex minus0.3ex
\parsep0.3ex plus0.2ex minus0.1ex}}{\end{list}}

\swapnumbers
\newtheorem{theorem}{Theorem}[section]
\newtheorem{lemma}[theorem]{Lemma}
\newtheorem{thm}[theorem]{}
\newtheorem{proposition}[theorem]{Proposition}
\newtheorem{definition}[theorem]{Definition}
\newtheorem{corollary}[theorem]{Corollary}
\newtheorem{remark}[theorem]{Remark}
\numberwithin{equation}{section}

\newcommand{\nG}{G}
\newcommand{\nH}{H}
\newcommand{\nT}{T}

\newcommand{\LRa}{\Leftrightarrow}

\newcommand{\A}{\mathbb{A}}
\newcommand{\B}{\mathbb{B}}
\newcommand{\cA}{\mathcal{A}}
\newcommand{\C}{\mathcal{C}}
\newcommand{\mC}{\mathbb{C}}
\newcommand{\M}{\mathbb{M}}
\newcommand{\X}{\mathbb{X}}
\newcommand{\Y}{\mathbb{Y}}
\newcommand{\bV}{\mathbb{V}}
\newcommand{\II}{\mathbb{I}}

\newcommand{\bF}{\mathbf{F}}
\newcommand{\G}{\mathbf{G}}
\newcommand{\T}{\mathbf{T}}
\newcommand{\HH}{\mathbf{HH}}
\newcommand{\bT}{\mathbf{T}}
\newcommand{\bH}{\mathbf{H}}

\newcommand{\AmT}{\A_\T}

\newcommand{\vareps}{\varepsilon}

\newcommand{\ttau}{\tilde{\tau}}
\newcommand{\oF}{{\overline{F}}}
\newcommand{\oH}{{\overline{H}}}
\newcommand{\woH}{{\widehat{\overline{H}}}}
\newcommand{\uH}{{\underline{H}}}
\newcommand{\wuH}{{\widehat{\underline{H}}}}
\newcommand{\ot}{\otimes}
\newcommand{\Hom}{{\rm Hom}}

\newcommand{\Nat}{{\rm Nat}}
\newcommand{\Inj}{{\rm {\bf Inj}}}

\newcommand{\cC}{{\mathcal C}}
\newcommand{\End}{{\rm End}}

\newcommand{\z}{\cdot}
\newcommand{\q}{{\bar{q}}}

\newcommand{\oa}{\overline{\alpha}}
\newcommand{\ob}{\overline{\beta}}
\newcommand{\og}{\overline{\gamma}}
\newcommand{\oR}{{\overline{R}}}
\newcommand{\uR}{\underline{R}}
\newcommand{\oL}{\overline{L}}
\newcommand{\uL}{\underline{L}}

\newcommand{\bR}{{\bf R}}
\newcommand{\ve}{\varepsilon}
\newcommand{\lra}{\longrightarrow}
\newcommand{\id}{I}

\begin{document}

\title{Bimonads and Hopf monads on categories}

\author{Bachuki Mesablishvili, Tbilisi\\ and \\ Robert Wisbauer, D\"usseldorf}

\maketitle

\begin{abstract}
The purpose of this paper is to develop a theory of bimonads and Hopf monads on arbitrary categories  thus providing the possibility
to transfer the essentials of the theory
of Hopf algebras in vector spaces to more general settings.
There are several extensions of this theory to
{\em monoidal} categories which in a certain sense follow the classical trace.
Here we do not pose any conditions on our base category but we
do refer to the monoidal structure of the category of endofunctors
on any category $\A$
and by this we retain some of the combinatorial complexity which makes
the theory so interesting. As a basic tool we use
distributive laws between monads and comonads (entwinings) on $\A$:
we define a {\em bimonad} on $\A$
as an endofunctor $B$ which is a monad and a comonad
with an entwining $\lambda:BB\to BB$ satisfying certain conditions.
This $\lambda$ is also employed to define the category $\A^B_B$ of
(mixed) $B$-bimodules.
In the classical situation, an entwining $\lambda$ is derived from the
twist map for vector spaces. Here this need not be the case
but there may exist special distributive laws $\tau:BB\to BB$
 satisfying the Yang-Baxter equation ({\em local prebraidings})
which induce an entwining $\lambda$ and lead to an
 extension of the theory of {\em braided Hopf algebras}.

An antipode is defined as a natural transformation $S:B\to B$ with
special properties and for categories $\A$ with limits or colimits
and bimonads $B$ preserving them, the existence of an antipode is
equivalent to $B$ inducing an equivalence between $\A$ and
the category $\A^B_B$ of $B$-bimodules. This is a general form of the
{\em Fundamental Theorem} of Hopf algebras.

Finally we observe a nice symmetry:
If $B$ is an endofunctor with a right adjoint $R$, then $B$
is a (Hopf) bimonad if and only if $R$ is a (Hopf) bimonad.
Thus a $k$-vector space $H$ is a Hopf algebra if and only if
$\Hom_k(H,-)$ is a Hopf bimonad. This provides a rich source for
Hopf monads not defined by tensor products and generalises the
well-known fact that a finite dimensional $k$-vector space $H$
is a Hopf algebra if and only if
 its dual $H^*=\Hom_k(H,k)$ is a Hopf algebra.
Moreover, we obtain that any set $G$ is a group if and only
if the functor $\rm{Map}(G,-)$ is a Hopf monad on the category of sets.
\end{abstract}

\tableofcontents

\section{Introduction}

The theory of algebras (monads) as well as of coalgebras (comonads) is well understood in
various fields of mathematis as algebra (e.g.\ \cite{BW}), universal algebra (e.g.\ \cite{Gumm}), logic or operational semantics (e.g.\ \cite{TuPl}), theoretical computer science (e.g.\ \cite{PoWa}).
The relationship between monads and comonads is controlled by
{\em distributive laws}
introduced in the seventies by Beck (see \cite{Beck}).
In algebra one of the fundamental notions emerging in this context are the Hopf algebras.
The definition is making heavy use of the tensor product and thus generalisations of this theory were mainly considered for {\em monoidal categories}.
They allow readily to transfer formalisms from
the category of vector spaces to the more general settings
(e.g. Bespalov and Brabant \cite{BeDr} and \cite{M}).

A Hopf algebra is an algebra as well as a coalgebra. Thus one way of
generalisation is to consider distinct algebras and coalgebras and some relationship
between them. This leads to the theory of {\em entwining structures} and
 {\em corings} over associative
rings (e.g.\ \cite{BW}) and one may ask how to formulate this in more
general categories.
The definition of {\em bimonads} on a monoidal category
as monads whose functor part is comonoidal
 by Brugui\`eres and Virelizier  in \cite[2.3]{BruVir}
may be seen as going in this direction.
Such functors are called {\em Hopf monads} in Moerdijk \cite{Moer} and {\em opmonoidal monads} in McCrudden \cite[Example 2.5]{McC}.
In \ref{opmon} we give more details of this notion.

Another extension of the theory of corings are
the {\em generalised bialgebras} in Loday in \cite{Lod}.
 These are Schur functors (on vector spaces) with a monad structure (operads)
and a specified coalgebra structure satisfying certain
compatibility conditions \cite[2.2.1]{Lod}.
While in \cite{Lod} use is made of the canonical twist map,
it is stressed in \cite{BruVir} that the theory is built up without
reference to any braiding. More comments on these constructions are given in
\ref{loday}.

The purpose of the present paper is to formulate the essentials
of the classical theory of Hopf algebras for any (not necessarily monoidal)
 category, thus making it accessible to a wide field of applications.
We also employ the fact
that the category of endofunctors (with the Godement product
as composition) always has a tensor product given by composition of natural transformations but no tensor product is required for the base category.

Compatibility between monads and comonads are
 formulated as distributive laws whose properties are recalled
in Section 2. In Section 3, general categorical notions
are presented and {\em Galois functors} are defined and investigated,
in particular equivalences induced for related categories
(relative injectives).

As suggested in \cite[5.13]{W-alg}, we define a {\em bimonad}
$\bH=(H,m,e,\delta,\vareps)$ on any category $\A$ as an endofunctor $H$
with a monad and a comonad structure satisfying compatibility conditions
(entwining) (see \ref{D.1.5}).
The latter do not refer to any braiding but in special cases they can
be derived from a {\em local prebraiding} $\tau:HH\to HH$ (see \ref{P.1.1}).
In this case the bimonad shows the characteristics of {\em braided bialgebras} (Section 6).

Related to a bimonad $H$ there is the (Eilenberg-Moore) category
$\A^\nH_\nH$ of bimodules with
a comparison functor $K_H: \A\to \A^\nH_\nH$.
An {\em antipode} is defined as a natural transformation $S:H\to H$ satisfying
$m \cdot  SH  \cdot \delta = e \cdot \varepsilon = m\cdot HS \cdot \delta$.
It exists if and only if the natural transformation
$\gamma := Hm\cdot \delta H:HH\to HH$ is an isomorphism.
If the category $\A$ is Cauchy complete
and $H$ preserves limits or colimits, the existence
of an antipode is equivalent to the comparison functor being an equivalence
(see \ref{T.1.14}). This is a general form of the Fundamental Theorem for
Hopf algebras. Any generalisation of Hopf algebras should offer an extension
of this important result.

Of course, bialgebras and Hopf algebras over commutative rings $R$ provide the
prototypes  for this theory: on $R$-Mod, the category of $R$-modules,
 one considers the endofunctor $B\ot_R -: R${\rm -Mod}$\to R${\rm -Mod}
where $B$ is an $R$-module with algebra and coalgebra structures, and
an entwining derived from the twist map (braiding) $M\ot_R N\to N\ot_R M$
(e.g.\ \cite[Section 8]{BBW}).

More generally, for a comonad $H$, the entwining $\lambda:HH\to HH$
may be derived from a {\em local prebraiding}  $\tau:HH\to HH$ (see \ref{Yang})
and then results similar to those known for braided Hopf algebras are obtained.
In particular, the composition $HH$ is again a bimonad  (see \ref{HH-bim})
and, if $\tau^2=1$, an {\em opposite bimonad} can be defined (see \ref{YB2}).

In case a bimonad $H$ on $\A$ has a right (or left) adjoint endofunctor
$R$, then $R$ is again a bimonad and has an antipode (or local prebraiding)
if and only if so does $H$ (see \ref{adj.bim}). In particular, for $R$-modules $B$,
the functor $\Hom_R(B,-)$ is right adjoint to $B\ot_R-$ and hence
$B$ is a Hopf algebra if and only if $\Hom_R(B,-)$ is a Hopf monad.
This provides a rich source for examples of Hopf monads not defined by a
tensor product and extends a symmetry principle known for finite dimensional
Hopf algebras (see \ref{Hopf-dual}). We close with the observation that
a set $G$ is a group if and only if the endofunctor $\rm{Map}(G,-)$
is a Hopf monad on the catgeory of sets (\ref{Hopf-group}).

Note that the pattern of our definition of bimonads resembles the definition
of {\em Frobenius monads} on any category by Street in \cite{Street}.
Those are monads ${\bf T}=(T,\mu,\eta)$
 with natural transformations $\ve:T\to \id$
and $\rho:T\to TT$, subject to suitable conditions,
 which induce a comonad structure $\delta= T\mu\cdot \rho T:T\to TT$ and
  product and coproduct on $T$ satisfy the compatibility condition
  $T\mu\cdot \delta T = \delta\cdot \mu = \mu T\cdot T\delta$.

\section{Distributive laws}

 Distributive laws between endofunctors were studied by Beck \cite{Beck},
Barr \cite{Barr} and others in the seventies of the last century.
They are a fundamental tool
for us and we recall some facts needed in the sequel. For more details and references we refer to \cite{W-alg}.

\begin{thm}\label{mon-com}{\bf Entwining  from monad  to comonad.} \em
Let $\T=(T,m,e)$ be a monad and $\G=(G,\delta,\vareps)$ a comonad on a category $\A$.
A natural transformation  $\lambda: TG \to GT$ is called a {\em mixed distributive law}
or {\em entwining} from the monad $\T$ to the comonad $\G$ if the diagrams

$$
\xymatrix{ & G \ar[dl]_{e G} \ar[rd]^{G e}&  \\
T G \ar[rr]_{\lambda}&& G T ,}  \qquad
\xymatrix{TG \ar[dr]_{T\varepsilon}\ar[rr]^\lambda  & & GT \ar[dl]^{\varepsilon T} \\
                           & T &} $$

$$
\xymatrix{ T G \ar[d]_{\lambda} \ar[r]^{T\delta} & T GG
\ar[r]^{\lambda G} & G T G \ar[d]^{G \lambda} &\text{and}& TT G
\ar[d]_{m G} \ar[r]^{T \lambda} & T G T \ar[r]^{\lambda T} & G T
T \ar[d]^{Gm}\\
G T \ar[rr]_{\delta T}&& G G T && TG \ar[rr]_{\lambda} && G T}
$$ are commutative.

It is shown in \cite{W} that for an arbitrary mixed distributive
law $\lambda: TG \to GT$
from a monad $\T$ to  a comonad $\G$, the triple $\widehat{\textbf{G}}=(\widehat{G},\widehat{\delta}, \widehat{\varepsilon})$,
 is a comonad on the category $\AmT$ of $\T$-modules (also called $\T$-algebras),
where  for any object $(a,h_a)$ of $\AmT$,
\begin{itemize}
\item $\widehat{G}(a, h_a)=(G(a), G(h_a) \cdot \lambda_a)$;
\item $(\widehat{\delta})_{(a, h_a)}=\delta_a$, and
\item $(\widehat{\varepsilon})_{(a, h_a)}=\varepsilon_a .$
\end{itemize}
 $\widehat{\textbf{G}}$ is called the {\em lifting of} $\textbf{G}$
corresponding to the mixed distributive law $\lambda$.

Furthermore, the triple $\widehat{\textbf{T}}=(\widehat{T}, \widehat{m}, \widehat{e})$
 is a monad on the category $\A^{\textbf{G}}$ of $\G$-comodules,
 where for any object $(a, \theta_a)$ of the category $\A^{\textbf{G}}$,

\begin{bull}
\item $\widehat{T}(a, \theta_a)=(T(a), \lambda_a \cdot T(\theta_a))$;
\item $(\widehat{m})_{(a, \theta_a)}=m_a$, and
\item $(\widehat{e})_{(a, \theta_a)}=e_a$.
\end{bull}

 This monad is called the {\em lifting of} $\textbf{T}$ corresponding to the mixed
 distributive law $\lambda$.
One has an isomorphism of categories $$(\A^{\nG})_{\widehat{\nT}}
\simeq (\A _{\nT})^{\widehat{\nG}},$$
and we write $\A^{\nG}_{\nT}(\lambda)$ for this category. An object of
$\A^{\nG}_{\nT}(\lambda)$ is a triple $(a, h_a, \theta_a)$, where
$(a, h_a) \in \A_{\nT}$ and $(a,\theta_a) \in \A^{\nG}$
with commuting diagram

\begin{equation}\label{D.1.1}
\xymatrix{
T(a) \ar[r]^-{h_a} \ar[d]_-{T(\theta_a)}& a \ar[r]^-{\theta_a}& G(a) \\
TG(a) \ar[rr]_-{\lambda_a}&& GT(a). \ar[u]_-{G(h_a)}}
\end{equation}
\end{thm}

We consider two examples of entwinings which may (also) be considered as
generalisations of Hopf algebras. They are different from our
approach and we will not refer to them later on.

 \begin{thm}\label{opmon} {\bf Opmonoidal functors.} \em
 Let $(\bV, \otimes, \II)$ be a strict monoidal category.
Following McCrudden \cite[Example 2.5]{McC}, one may call
a monad $(T,\mu,\eta)$ on $\bV$ {\em opmonoidal} if there exist
  morphisms
\begin{center}
$\theta: T(\II)\to \II$\; and  \; $\chi_{X,Y}: T(X\ot Y)\to T(X)\ot T(Y)$,
\end{center}
the latter natural in $X,Y\in \bV$,
which are compatible with the tensor structure of $\bV$ and the monad structure
of $T$.

Such functors can also be characterised by the condition that the
tensor product of $\bV$ can be lifted to the category of $T$-modules
(e.g. \cite[3.4]{W-alg}).
They were introduced and named {\em Hopf monads} by Moerdijk in
\cite[Definition 1.1]{Moer} and called {\em bimonads} by
Brugui\`eres and Virelizier in \cite[2.3]{BruVir}.
It is mentioned in \cite[Example 2.8]{BruVir} that Szlach\'anyi's
bialgebroids in \cite{SzlA} may be interpreted in terms of such "bimonads".
It is preferable to use the
terminology from  \cite{McC} since these functors are neither
bimonads nor Hopf monads in a strict sense but rather an entwining
(as in \ref{mon-com}) between the monad $T$ and the comonad
$T(\II)\ot-$ on $\bV$:

Indeed, the compatibility conditions required in the definitions
induce a coproduct $\chi_{\II,\II}: T(\II)\to T(\II)\ot T(\II)$ with
counit $\theta: T(\II)\to \II$. Moreover, the relation between
$\chi$ and $\mu$ (e.g. (15) in \cite[2.3]{BruVir}) lead to the
commutative diagram
(using $X\ot \II=X$)
$$\xymatrix{
 TT(X) \ar[rr]^\mu \ar[d]_{T(\chi_{\II,X})} & &
        T(X) \ar[rr]^{\chi_{\II,X}} & & T(\II)\ot T(X) \\
T(T(\II)\ot T(X)) \ar[rr]^{\chi_{T(\II),T(X)}}& & TT(\II)\ot TT(X)
  \ar[rr]^{\mu_\II \ot TT(X)}& &
                 T(\II)\ot TT(X) \ar[u]_{T(\II)\ot \mu_X}
}$$
This shows that $T(X)$ is a mixed $(T,T(\II)\ot-)$-bimodule for the
entwining map
 $$\lambda=  (\mu_\II \ot T(-)) \circ\chi_{T(\II),-}:
    T(T(\II)\ot -) \to T(\II)\ot T(-). $$
The {\em antipode} of a classical Hopf algebra $H$ is defined as a special
endomorphism of $H$. Since opmonoidal monads $T$ relate two distinct functors
it is not surprising that the notion of an antipode can not
be transferred easily to this situation and the attempt to do so leads
to an "apparently complicated definition" in \cite[3.3 and Remark 3.5]{BruVir}.
Hereby the base category $\C$ is required to be {\em autonomous}.
\end{thm}

\begin{thm}\label{loday} {\bf Generalized bialgebras and Hopf operads.} \em
The {\em generalised bialgebras} over fields as defined in
Loday \cite[Section 2.1]{Lod}
are similar to the mixed bimodules (see \ref{mon-com}):
they are vector spaces which are modules over some {\em operad ${\cA}$} (Schur functors with multiplication and unit) and comodules over some
coalgebras $\cC^c$,  which are linear duals of some operad $\C$.
Similar to the opmonoidal monads the coalgebraic structure is
based on the tensor product (of vector spaces).
 The Hypothesis (H0) in \cite{BruVir} resembles the role of the
entwining $\lambda$ in \ref{D.1.1}.
The Hypothesis (H1) requires that the free ${\cA}$-algebra is a $(\cC^C,{\cA})$-bialgebra:
this is similar to the condition on an $A$-coring $C$, $A$ an associative algebra, to have a $C$-comodule structure (equivalently the existence of a group-like element, e.g. \cite[28.2]{BW}). The condition (H2iso) plays the role of
the canonical isomorphism defining {\em Galois corings} and
the {\em Galois Coring Structure Theorem} \cite[28.19]{BW}
may be compared with the
{\em Rigidity Theorem} \cite[2.3.7]{Lod}. The latter
can be considered as a generalisation of the Hopf-Borel Theorem (see \cite[4.1.8]{Lod}) and of the Cartier-Milnor-Moore Theorem
(see \cite[4.1.3]{Lod}).
In \cite[3.2]{Lod}, {\em Hopf operads} are defined in the sense of Moerdijk \cite{Moer} and thus the coalgebraic part is dependent on the tensor
product.
This is only a sketch of the similarities between Loday's setting and
our approach here. It will be interesting to work out the
relationship in more detail.
\end{thm}

Similar to \ref{mon-com} we will also need the notion of mixed distributive
laws from a comonad to a monad.

\begin{thm}\label{common} {\bf Entwining from comonad to monad.} \em
A natural transformation $\lambda : GT \to TG$ is a {\em mixed distributive law} from a comonad $\G$ to a monad $\T$, also called an {\em entwining} of $\G$ and $\T$, if the   diagrams

$$\xymatrix{
& G \ar[dl]_-{Ge} \ar[rd]^-{eG}& &GT \ar[dr]_-{\varepsilon T} \ar[rr]^-{\lambda}&& TG
 \ar[ld]^-{T \varepsilon}\\
GT \ar[rr]_-{\lambda}&& TG\,,  && T&}
$$

$$
\xymatrix{ GTT \ar[d]_{Gm} \ar[r]^{\lambda T} & T G T
\ar[r]^{T \lambda } & T T G \ar[d]^{mG}&
 GGT
\ar[r]^{G \lambda} &  G T G \ar[r]^{\lambda G} &
TGG \\
G T \ar[rr]_{\lambda}&& TG, &  GT \ar[u]^{\delta T}\ar[rr]_{\lambda} && TG \ar[u]_{T\delta}}
$$ are commutive.
\end{thm}

For convenience we recall the distributive laws between two monads
and between two comonads   (e.g. \cite{Beck}, \cite{Barr}, \cite[4.4 and 4.9]{W-alg}).

 \begin{thm} \label{lift.mon.mon}  {\bf Monad distributive.} \em
Let $\bF=(F,m,e)$ and $\T=(T,m',e')$ be monads
  on the category $\A$.
A natural transformation  $\lambda: FT\to TF$
is said to be {\em monad distributive} if it induces the commutative diagrams

 $$
    \xymatrix{ & T\ar[dl]_{e_T} \ar[dr]^{Te}& \\
       FT \ar[rr]^\lambda  & & TF,} \qquad
 \xymatrix{ & F\ar[dl]_{Fe'} \ar[dr]^{e'_F} & \\
   FT \ar[rr]^\lambda & & TF.}
  $$
  $$
  \xymatrix{
   FFT\ar[rr]^{m_T} \ar[d]_{F\lambda} & & FT \ar[d]^\lambda \\
   FTF \ar[r]^{\lambda_F} & TFF \ar[r]^{Tm} & TF , } \qquad
      \xymatrix{
   FTT\ar[rr] ^{Fm'} \ar[d]_{\lambda_T} & & FT \ar[d]^\lambda \\
   TFT \ar[r]^{T\lambda} & TTF \ar[r]^{m'_F} & TF .}
      $$
     In this case
   $\lambda:FT\to TF $ induces a canonical monad structure on $TF$.
\end{thm}

  \begin{thm} \label{lift.comon.comon}  {\bf Comonad  distributive.} \em
Let $\G=(G,\delta,\vareps)$ and  $\T=(T,\delta',\vareps')$ be comonads
 on the category $\A$.
A natural transformation $\varphi: TG\to GT$ is said to be {\em comonad distributive}
if it induces the commutative diagrams

$$
 \xymatrix{ TG\ar[dr]_{T\vareps} \ar[rr]^{\varphi} && GT \ar[dl]^{\vareps_T} \\
         & T &,} \qquad
 \xymatrix{ TG\ar[dr]_{\vareps'_G} \ar[rr]^{\varphi} & & GT \ar[dl]^{G\vareps'} \\
    & G &, } $$
$$
\xymatrix{
   TG\ar[r] ^{T\delta} \ar[d]_{\varphi}  & TGG \ar[r]^{\varphi_G} & GTG \ar[d]^{G\varphi} \\
   GT \ar[rr]^{\delta_T} & &  GGT ,} \qquad
     \xymatrix{
   TG\ar[r] ^{\delta'_G} \ar[d]_{\varphi}  & TTG \ar[r]^{T\varphi} & TGT \ar[d]^{\varphi_T} \\
   GT \ar[rr]^{G\delta'} & &  GTT    .  }
 $$

 In this case
    $\varphi\!:\! TG\to GT$ induces a canonical comonad structure on $TG$.
\end{thm}

 \section{Actions on functors and Galois functors}

The language of modules over rings can also be used to describe actions
of monads on functors. Doing this we define Galois functors and
to characterise those we investigate the
relationships between categories of relative injective objects.

\begin{thm}\label{T-act}{\bf $\T$-actions on functors.} \em
Let $\A$ and $\B$ be categories.
Given a monad $\T=(T,m,e)$ on $\A$ and any functor $L : \A \to \B$, we say that $L$ is a
{\em (right) $\textbf{T}$-module} if there exists a natural transformation
$\alpha_L : LT \to L$ such that the diagrams
\begin{equation}\label{E.1.21}
\xymatrix{
L \ar@{=}[dr] \ar[r]^-{Le}&LT \ar[d]^-{\alpha_L}\\
& L ,} \qquad
\xymatrix{
LTT \ar[r]^-{L m} \ar[d]_-{\alpha_L T}& LT \ar[d]^-{\alpha_L}\\
LT \ar[r]_-{\alpha_L}& L}
\end{equation} commute. It is easy to see that $(T,m)$ and $(TT, Tm)$
both are $\textbf{T}$-modules.

Similarly, given a comonad $\G=(G,\delta,\vareps)$ on $\A$, a functor $K: \B \to
\A$ is a {\em left $\textbf{G}$-comodule} if there exists a
natural transformation $\beta_K : K \to GK$ for which the diagrams

$$
\xymatrix{
K \ar@{=}[dr] \ar[r]^-{\beta_K}& GK \ar[d]^-{\varepsilon K}\\
& K ,} \qquad
\xymatrix{
K \ar[r]^-{\beta_K} \ar[d]_-{\beta_K}& GK \ar[d]^-{\delta K}\\
GK \ar[r]_-{G \beta_K}& GGK}
$$ commute.

Given two $\textbf{T}$-modules $(L, \alpha_L)$, $(L', \alpha_{L'})$, a natural transformation
$g: L \to L'$ is called {\em $\textbf{T}$-linear} if the diagram
\begin{equation}\label{E.1.23}
\xymatrix{
LT \ar[r]^-{gT} \ar[d]_-{\alpha_L}& L'T \ar[d]^-{\alpha_{L'}}\\
L \ar[r]_{g} & L'}
\end{equation} commutes.
\end{thm}

\begin{lemma}\label{L.1.9}
Let $(L, \alpha_L)$ be a $\textbf{T}$-module. If $f, f': TT \to L$ are
$\textbf{T}$-linear morphisms from the $\textbf{T}$-module $(TT, Tm)$ to
the $\textbf{T}$-module $(L, \alpha_L)$ such that $f \cdot Te=f' \cdot Te$, then $f=f'$.
\end{lemma}

\begin{proof}
Since $f \cdot Te=f' \cdot Te$, we have $ \alpha_{L} \cdot fT \cdot TeT=\alpha_{L}\cdot f'T \cdot TeT.$ Moreover, since $f$ and $f'$ are both $\textbf{T}$-linear, we have the commutative
diagrams
$$\xymatrix{
TTT \ar[d]_{Tm}\ar[r]^{fT}& LT \ar[d]^{\alpha_L}\\
TT \ar[r]^{f}&L ,
}\qquad
\xymatrix{
TTT \ar[d]_{Tm}\ar[r]^{f'T}  & LT \ar[d]^{\alpha_L}\\
TT \ar[r]^{f'} &L.
}
$$
Thus $\alpha_{L} \cdot fT=f \cdot Tm$ and $\alpha_{L} \cdot f'T=f' \cdot Tm,$ and we have $f \cdot Tm \cdot TeT=f' \cdot Tm \cdot TeT$. It follows - since $Tm \cdot TeT=1$ - that $f=f'$.
\end{proof}

\begin{thm}\label{left-G}{\bf Left $G$-comodule functors.} \em
Let $\G$ be a comonad on a category $\A$, let $U^{\nG}:
\A^{\nG}  \to \A$ be the forgetful functor and write
 $\phi^G:\A\to \A^\nG$ for the cofree $\nG$-comodule functor.
 Fix a functor $F: \B \to \A$, and consider a functor
$\overline{F} :\B\to \A^{\nG}$ making the diagram
\begin{equation}\label{D.3.1}
\xymatrix{\B \ar[rr]^-{\overline{F}} \ar[dr]_{F} && \A^{G}
\ar[dl]^{U^{G}}\\
& \A & } \end{equation}
commutative. Then $\overline{F} (b)=(F(b),
\alpha_{F(b)})$ for some $\alpha_{F(b)} : F(b) \to GF(b)$.
Consider the natural transformation
\begin{equation} \bar{\alpha}_F : F \to GF,
\end{equation}
whose $b$-component is $\alpha_{F(b)}$. It should
be pointed out that $\bar{\alpha}_F $ makes $F$ a left
$\textbf{G}$-comodule, and it is easy to see that there is a one
to one correspondence between functors
$\overline{F} :\B \to\A^{G}$
making the diagram (\ref{D.3.1}) commute and natural
transformations $\bar{\alpha}_F : F \to GF$ making $F$ a left
$\textbf{G}$-comodule.
\end{thm}

The following is an immediate consequence of (the dual of)
\cite[Propositions II,1.1 and II,1.4]{D}:

\begin{theorem}\label{T.1}
Suppose that $F$ has a right adjoint $R : \A \to \B$ with unit $
\eta: 1 \to RF$ and counit $\varepsilon : FR \to 1$. Then the
composite
$$\xymatrix{t_{\overline{F}}:
FR \ar[r]^-{\bar{\alpha}_{F} R} & GFR \ar[r]^{G \varepsilon}
&G.}$$ is a morphism from the comonad $\textbf{G}'=(FR,F\eta R,\ve)$
generated by the adjunction $\eta,
\varepsilon : F \dashv R : \A \to \B$ to the comonad $\textbf{G}$.
Moreover, the assignment
$$\overline{F} \longrightarrow t_{\overline{F}}$$ yields a one to one
correspondence between functors $\overline{F} : \B \to \A^\nG$ making the
diagram (\ref{D.3.1}) commutative and morphisms of comonads $t_{\overline{F}} :
\textbf{G}' \to \textbf{G}$.
\end{theorem}

\begin{definition} \em We say that a left $\textbf{G}$-comodule $F: \B \to
\A$ with a right adjoint $R: \B \to \A$ is {\em $\textbf{G}$-Galois} if
the corresponding morphism $t_{\overline{F}}: FR \to \textbf{G}$ of
comonads on $\A$ is an isomorphism.
\end{definition}

As an example, consider an $A$-coring $\cC$, $A$ an associative ring,
and any right $\cC$-comodule $P$ with $S=\End^\cC(P)$.
Then there is a natural transformation
$$ \tilde \mu: \Hom_A(P,-)\ot_S P \to -\ot_A \cC$$
and $P$ is called a {\em Galois comodule} provided
$ \tilde \mu_X$ is an isomorphism for any right $A$-module $X$,
that is, the functor
$-\ot_S P: \M_S\to \M^\cC$ is a $-\ot_A\cC$-Galois comodule
(see \cite[Definiton 4.1]{WisGal}).
\medskip

\begin{thm}\label{right.adj}{\bf Right adjoint functor of $\oF$.} \em
When the category $\B$ has equalisers, the functor $\overline{F}$
has a right adjoint, which can be described as follows: Writing
$\beta_R$ for the composite
$$\xymatrix{R \ar[r]^{\eta R\quad}& RFR \ar[r]^{R t_{\overline{F}}}& RG,}$$
it is not hard to see that the equaliser $(\overline{R},
\overline{e})$ of the following diagram
$$\xymatrix{ RU^G
\ar@{->}@<0.5ex>[rr]^-{RU^G \eta_G} \ar@ {->}@<-0.5ex>
[rr]_-{\beta_R U^G}&& RGU^G=RU^G\phi^GU^G,}$$ where $\eta_G : 1 \to
\phi^GU^G$ is the unit of the adjunction $U^G \dashv \phi^G,$ is right
adjoint to $\overline{F}$.
\end{thm}

\begin{thm}\label{MoB.0.1}{\bf Adjoints and monads.} \em
For  categories $\A$, $\B$,
let $L:\A\to \B$ be a functor with right adjoint $R:\B\to \A$.
Let $\bT=(T,m, e)$ be a monad on $\A$ and
suppose there exists a functor 
$\oR:\B \to \A_T$ yielding the commutative diagram
$$
\xymatrix{ \B \ar[rr]^\oR \ar[rd]_R & & \A_T \ar[ld]^{U_T}\\
     & \A.}$$  
Then  $\oR(b)=(R(b), \beta_{b})$ for some
$\beta_{b}: TR(b) \to R(b)$ and the collection 
$\{\beta_b,\, b \in \B\}$ 
constitutes a natural transformation $\beta_\oR: TR \to R$. 
It
is proved in \cite{D} that the natural transformation
$$t_\oR: 
\xymatrix{T \ar[r]^-{T \eta}& TRL \ar[r]^-{\beta L} & RL}$$ 
is a
morphism of monads.  By the dual of \cite[Theorem 4.4]{M},
we obtain:
\smallskip

{\em The functor $\oR$ is an equivalence of categories
iff the functor $R$ is monadic and $t_\oR$ is an isomorphism of
monads.}
\end{thm}

In view of the characterisation of Galois functors we have a closer
look at some related classes of relative injective objects.

Let $F : \B \to \A$ be any functor. Recall (from \cite{H}) that an object
$b \in \B$ is said to be $F$-\textit{injective} if for any diagram
in $\B$,
$$\xymatrix{b_1 \ar[d]_g \ar[r]^f & b_2 \ar@{.>}[dl]^h\\
b&}$$
with $F(f)$ a split monomorphism in $\A$, there exists a
morphism $h: b_2 \to b$ such that $hf=g$. We write
$\Inj (F, \B)$ for the full subcategory of $\B$ with
objects all $F$-injectives.

The following result from \cite{S} will be needed.

\begin{proposition}\label{P.3.3}
Let $\eta, \varepsilon: F \dashv R : \A \to \B$ be an adjunction.
For any object $b \in \B$, the following assertions are
equivalent:

\begin{blist}
\item  $b$ is $F$-injective;
\item $b$ is a coretract for some $R(a)$, with $a \in \A$;
\item the $b$-component $\eta_b : b \to RF(b)$ of $\eta$ is
      a split monomorphism.
\end{blist}
\end{proposition}

\begin{remark}\label{R.3.4} \em
For any $a \in ¸\A$, $R(\varepsilon_a)\cdot \eta_{R(a)}=1$ by one
of the triangular identities for the adjunction $F \dashv R$.
Thus, $R(a) \in \Inj(F, \B)$ for all $a \in \A$.
Moreover, since the composite of coretracts is again a coretract,
it follows from (b) that $\Inj(F, \B)$ is closed
under coretracts.
\end{remark}

\begin{thm}\label{P.3.5}{\bf Functor between injectives.} \em
Let $K_{G'}: \B \to \A^{\nG'}$ be the comparison functor
(notation as in \ref{T.1}).
 If $b\in \B$ is $F$-injective, then $K_{G'}(b)=(F(b), F(\eta_b))$ is
$U_{G'}$-injective, since by the fact
that $\eta_b$ is a split monomorphism in $\B$,
 $(\eta_{G'})_{\phi^{G'}(b)}=F(\eta_b)$
is a split monomorphism in $\A^{\nG'}$ ($G'$ as in \ref{T.1}).
Thus the functor $K_{G'}: \B \to \A_{G'}$ yields a functor
$$ \Inj (K_{G'}):  \Inj(\bF,\B) \to \Inj(\phi^{G'},\A^{\nG'}).$$
  When $\B$ has equalisers, this functor 
is an equivalence of categories (see \cite{S}).
\end{thm}

We shall henceforth assume that $\B$ has equalisers.

\begin{proposition}\label{P.3.6} The functor $\overline{R}: \A^\nG \to \B$
restricts to a functor $$\overline{R}': \Inj(U^{G},
\A^{\nG}) \to \Inj(F, \B).$$
\end{proposition}
\begin{proof} Let $(a, \theta_a)$ be an arbitrary object of
$\Inj (U^{G}, \A^{\nG})$. Then, by Proposition \ref{P.3.3},
there exists an object $a_0 \in \A$ such that $(a, \theta_a)$ is a
coretraction of $\phi^G(a_0)=(G(a_0), \delta_{a_0})$ in $\A^\nG$, i.e.,
there exist morphisms
\begin{center}
$f: (a, \theta_a) \to (G(a_0),
\delta_{a_0})$ and $g: (G(a_0), \delta_{a_0}) \to (a, \theta_a)$
\end{center}
in $\A^\nG$ with $gf=1$. Since $f$ and $g$ are morphisms in $\A^\nG$,
the diagram
$$\xymatrix{
G(a_0) \ar@{->}@<+1.5ex>[d]^-{g} \ar[r]^{(\delta_G)_{a_0}} &
GG(a_0)\ar@{->}@<+1.5ex>[d]^-{G(g)}\\
a \ar@{->}@<-0.5ex>[u]^-{f} \ar[r]_{\theta_a}&
G(a)\ar@{->}@<-0.5ex>[u]^-{G(f)} }
$$commutes. By naturality of $\beta_R$, the  diagram
$$\xymatrix{
RG(a_0) \ar@{->}@<+1.5ex>[d]^-{R(g)} \ar[rr]^{(\beta_R)_{G(a_0)}}
&&RGG(a_0)
\ar@{->}@<+1.5ex>[d]^-{RG(g)}\\
R(a) \ar@{->}@<-0.5ex>[u]^-{R(f)} \ar[rr]_{(\beta_R)_a}&&
RG(a)\ar@{->}@<-0.5ex>[u]^-{RG(f)} }
$$ also commutes. Consider now the following commutative diagram
\begin{equation}\label{D.3.3}\xymatrix{R(a_0) \ar@{.>}@<+1.5ex>[dd] \ar[r]^{\beta_{a_0}}&
RG(a_0) \ar@{->}@<+1.5ex>[dd]^-{R(g)}
\ar@{->}@<0.5ex>[rr]^{(\beta_R)_{G(a_0)}}
\ar@{->}@<-0.5ex>[rr]_{R((\delta_G)_{a_0})} &&RGG(a_0)\ar@{->}@<+1.5ex>[dd]^-{RG(g)}\\\\
\overline{R}(a, \theta_a)\ar@{.>}@<-0.5ex>[uu]
\ar[r]_{\overline{e}_{(a, \theta_a)}}&R(a)
\ar@{->}@<-0.5ex>[uu]^-{R(f)} \ar@{->}@<0.5ex>[rr]^{(\beta_R)_a}
\ar@{->}@<-0.5ex>[rr]_{R(\theta_a)}&&
RG(a)\ar@{->}@<-0.5ex>[uu]^-{RG(f)} .} \end{equation}
 It is not
hard to see that the top row of this diagram is a (split)
equaliser (see, \cite{G}), and since the bottom row is an
equaliser by the very definition of $\overline{e}$, it follows
from the commutativity of the diagram that $\overline{R}(a,
\theta_a)$ is a coretract of $R(a_0)$, and thus is an object of
$\Inj (F, \B)$ (see Remark \ref{R.3.4}). It means that the
functor $\overline{R}: \A^\nG \to \B$ can be restricted to a functor
$\overline{R}': \Inj (U^G, \A^\nG)\to
\Inj (F, \B)$.
\end{proof}

\begin{proposition}\label{P.3.7}
Suppose that for any $b \in \B$, $(t_{\overline{F}})_{F(b)}$ is an
isomorphism. Then the functor $\overline{F}: \B \to \A^\nG$ can be
restricted to a functor
$$\overline{F}': \Inj (F, \B)\to \Inj (U^G, \A^\nG).$$
\end{proposition}

\begin{proof} Let $\delta'$ denote the comultiplication in the
comonad $\G'$ (see \ref{T.1}), i.e., $\delta'= F \eta R.$ Then for any $b \in \B$,
$$\begin{array}{rl}
\overline{F}(RF(b))&=\A_{t_{\overline{F}}}(\phi^{G'}(UF(b)))
   =\A_{t_{\overline{F}}}(FRF(b), F \eta_{RF(b)})\\[+1mm]
& =\A_{t_{\overline{F}}}(G'F(b), \delta'_{F(b)})=
 (G'F(b),(t_{\overline{F}})_{G'F(b)}\cdot\delta'_{F(b)}).
\end{array}$$
Consider now the diagram

$$
\xymatrix{ G'F(b) \ar[rrr]^-{(t_{\overline{F}})_{F(b)}}
\ar[d]_{\delta'_{F(b)}}&&& GF(b) \ar[ddd]^{\delta_{F(b)}}\\
G'G'F(b)
\ar@{}[rrru]_{(1)}\ar[rrrdd]^{(t_{\overline{F}})_{F(b)}.(t_{\overline{F}})_{F(b)}}
\ar[dd]_{(t_{\overline{F}})_{G'F(b)}} &&& &\\\\
GG'F(b) \ar[rrr]_{G((t_{\overline{F}})_{F(b)})}&&& GGF(b)\, ,}
$$ in which the triangle commutes by the 
definition of the
composite $(t_{\overline{F}})_{F(b)}.(t_{\overline{F}})_{F(b)}$, while the
diagram (1) commutes since $t_{\overline{F}}$ is a morphism of
comonads.
The commutativity of the outer diagram shows that
$(t_{\overline{F}})_{F(b)}$ is a morphism from the $G$-coalgebra
$\overline{F}(RF(b))=(G'F(b),(t_{\overline{F}})_{G'F(b)}\cdot
\delta'_{F(b)})$ to the $G$-coalgebra $(GF(b), \delta_{F(b)})$.
Moreover, $(t_{\overline{F}})_{F(b)}$ is an isomorphism by our
assumption. Thus, for any $b \in \B$, $\overline{F}(RF(b))$ is
isomorphic to the $G$-coalgebra $(GF(b), \delta_{F(b)})$, which is
of course an object of the category $\Inj (U^G,
\A^\nG)$. Now, since any $b \in \Inj (F, \B)$ is a
coretract of $RF(b)$ (see Remark \ref{R.3.4}), and since any functor takes
coretracts to coretracts, it follows that, for any $b \in
\Inj (F, \B)$, $\overline{F}(b)$ is a coretract of
the $G$-coalgebra $(GF(b), \delta_{F(b)})\in
\Inj (U^G, \A^\nG)$, and thus is an object of the
category $\Inj (U^G, \A^\nG)$ again by Remark \ref{R.3.4}.
This completes the proof.
\end{proof}

The following technical observation is needed for the next proposition.

\begin{lemma}\label{L.3.8} Let $\iota, \kappa : W \dashv W' : \Y \to \X$ be an
adjunction of any categories. If $i: x' \to x$ and $j: x \to x'$ are morphisms in
$\X$ such that $ji=1$ and if $\iota_x$ is an isomorphism, then
$\iota_{x'}$ is also an isomorphism.
\end{lemma}

\begin{proof} Since $ji=1$, the  diagram
$$\xymatrix{ x'\ar[r]^{i}&x
\ar@{->}@<0.5ex>[r]^-{1} \ar@ {->}@<-0.5ex> [r]_-{ij}& x}$$ is a
split equaliser. Then the diagram $$\xymatrix{ W'W(x')
\ar[rr]^-{W'W(i)}&& W'W(x) \ar@{->}@<0.5ex>[rr]^-{1} \ar@
{->}@<-0.5ex> [rr]_-{W'W(ij)}&& W'W(x)}$$ is also a split
equaliser. Now considering the following commutative diagram

$$
\xymatrix{x' \ar@{.>}[d]_{\iota_{x'}}\ar[rr]^{i}&& x
\ar[d]_{\kappa_x} \ar@{->}@<0.5ex>[rr]^-{1} \ar@
{->}@<-0.5ex> [rr]_-{ij}&& x \ar[d]^{\kappa_x}\\
W'W(x') \ar[rr]_-{W'W(i)}&& W'W(x) \ar@{->}@<0.5ex>[rr]^-{1} \ar@
{->}@<-0.5ex> [rr]_-{W'W(ij)}&& W'W(x)}$$ and recalling that the
vertical two morphisms are both isomorphisms by assumption, we get
that the morphism $\iota_{x'}$ is also an isomorphism.
\end{proof}

\begin{proposition}\label{P.3.9} In the situation of Proposition \ref{P.3.7},
$\Inj (F,\B)$ is (isomorphic to) a coreflective subcategory of the category
$\Inj (U^G, \A^\nG)$.
\end{proposition}

\begin{proof} By Proposition \ref{P.3.6}, the functor
$\overline{R}$ restricts to a functor
$$\overline{R}':\Inj (U^G, \A^\nG) \to \Inj (F, \B),$$
while according to Proposition \ref{P.3.7},
the functor $\overline{F}$ restricts to a functor
$$\overline{F}':\Inj (F, \B) \to \Inj (U^G,\A^\nG).$$
Since
\begin{itemize}
\item $\overline{F}$ is a left adjoint to $\overline{R}$,
\item $\Inj (F, \B)$ is a full subcategory of $\B$, and
\item $\Inj (U^G, \A^\nG)$ is a full subcategory of $\A^\nG$,
\end{itemize}
the functor $\overline{F}'$ is left adjoint
to the functor $\overline{R}'$, and the unit $\overline{\eta}': 1
\to \overline{R}' \overline{F}'$ of the adjunction $\overline{F}'
\dashv \overline{R}'$ is the restriction of $\overline{\eta}:
\overline{F} \dashv \overline{R}$ to the subcategory
$\Inj (F, \B)$, while the counit
$\overline{\varepsilon}': \overline{F}' \overline{R}' \to 1$ of
this adjunction is the restriction of  $\overline{\varepsilon}:
\overline{F}\overline{R}\to 1$ to the subcategory
$\Inj (U^G, \A^\nG)$.

Next, since the top of the diagram \ref{D.3.3} is a (split) equaliser,
$\overline{R}(G(a_0), \delta_{a_0})\simeq R(a_0)$. In particular,
taking $(GF(b), \delta_{F(b)})$, we see that
$$RF(b) \simeq
\overline{R}(GF(b), \delta_{F(b)})= \overline{R} \,
\overline{F}(UF(b)).$$
 Thus, the $RF(b)$-component
$\overline{\eta}'_{RF(b)}$ of the unit $\overline{\eta}': 1 \to
\overline{R}' \overline{F}'$ of the adjunction $\overline{F}'
\dashv \overline{R}'$ is an isomorphism. It now follows from Lemma
\ref{L.3.8} - since any $b \in \Inj (F, \B)$ is a
coretraction of $RF(b)$ - that $\overline{\eta}'_b$ is an
isomorphism for all $b \in \Inj (F, \B)$ proving
that the unit $\overline{\eta}'$ of the adjunction $\overline{F}'
\dashv \overline{R}'$ is an isomorphism. Thus
$\Inj (F, \B)$ is (isomorphic to) a coreflective
subcategory of the category $\Inj (U^G, \A^\nG)$.
\end{proof}

\begin{corollary}\label{C.3.10} In the situation of Proposition \ref{P.3.7},
suppose that each component of the unit $\eta: 1 \to RF$ is a split
monomorphism. Then the category $\B$ is (isomorphic to) a
coreflective subcategory of $\Inj (U^G, \A^\nG)$.
\end{corollary}

\begin{proof} When each component of the unit $\eta: 1 \to RF$ is a split
monomorphism, it follows from Proposition \ref{P.3.3} that every $b \in
\B$ is $F$-injective; i.e. $\B=\Inj (F, \B)$. The
assertion now follows from Proposition \ref{P.3.9}.
\end{proof}

\begin{thm}\label{T.3.11}{\bf Characterisation of $\G$-Galois comodules.}
Assume $\B$ to admit equalisers,  let $\G$ be a comonad on $\A$, and
$F:\B\to \A$ a functor with right adjoint $R:\A\to \B$. If there
exists a functor $\overline{F} : \A \to \A^{\nG}$ with
$U^{\nG}\overline{F}=F$, then the following are equivalent:
\begin{blist}
\item  $F$ is $\G$-Galois, i.e. $t_{\overline{F}}: \textbf{G}\,'\to
        \textbf{G}$ is an isomorphism;
\item the following
composite is an isomorphism:
$$\xymatrix{ \overline{F} R \ar[r]^-{\eta_G
\overline{F}R} &\phi^G U^G \overline{F}R=\phi^G FR \ar[r]^-{\phi^G
\varepsilon}& \phi^G};$$
\item the functor
$\overline{F}: \B \to \A^{\nG}$ restricts to an equivalence of
categories $$\Inj (F, \B) \to
\Inj (U^{G}, \A^{\nG});$$
\item  for any $(a,\theta_a) \in \Inj (U^{G}, \A^{\nG})$, the $(a,
\theta_a)$-component $\overline{\varepsilon}_{(a, \theta_a)}$ of
the counit $\overline{\varepsilon}$ of the adjunction
$\overline{F} \dashv \overline{R}$, is an isomorphism;
\item
for any $a \in \A$,
$\overline{\varepsilon}_{\phi_{G}(a)}=\overline{\varepsilon}_{(G(a),
\delta_a)}$ is an isomorphism.
\end{blist}
\end{thm}

\begin{proof}
That (a) and (b) are equivalent is proved in \cite{D1}.
By the proof of \cite[Theorem of 2.6]{G},
for any $a \in \A$,
$\overline{\varepsilon}_{\phi^G (a)}=\overline{\varepsilon}_{(G(a),
\delta_a)}=(t_{\overline{F}})_a$, thus (a) and (e) are equivalent.

By Remark \ref{R.3.4}, (d) implies (e).

Since $\B$ admits equalisers by our assumption on $\B$, it follows
from Proposition \ref{P.3.5} that the functor
$\Inj (K_{G'})$ is an equivalence of categories.
Now, if $t_{\overline{F}}: \textbf{G} ' \to \textbf{G}$ is an
isomorphism of comonads, then the functor $\A_{t_{\overline{F}}}$
is an isomorphism of categories, and thus $\overline{F}$ is
isomorphic to the comparison functor $K_{G'}$. It now follows from
Proposition \ref{P.3.5} that $\overline{F}$ restricts to the functor
$\Inj (F, \B) \to \Inj (U^{G},\A^{\nG})$
 which is an equivalence of categories.
Thus (a) $\Rightarrow$ (c).

If the functor $\overline{F}: \B \to \A^\nG$ restricts to a functor
$$\overline{F}':\Inj (F, \B) \to \Inj (U^{G},\A^{\nG}),$$
then one can prove as in the proof of Proposition 3.9
that $\overline{F}'$ is left adjoint to $\overline{R}'$ and that
the counit $\overline{\varepsilon}': \overline{F}' \,
\overline{R}' \to 1$ of this adjunction is the restriction of the
counit $\overline{\varepsilon}: \overline{F}\, \overline{R} \to 1$
of the adjunction $\overline{F} \dashv \overline{R}$ to the
subcategory $\Inj (U^{G}, \A^{\nG})$. Now, if
$\overline{F}'$ is an equivalence of categories, then
$\overline{\varepsilon}'$ is an isomorphism. Thus, for any $(a,
\theta_a) \in \Inj (U^{G}, \A^{G})$,
$\overline{\varepsilon}'_{(a, \theta_a)}$ is an isomorphism
proving that (c)$\Rightarrow$(d).
\end{proof}

\section{Bimonads}

The following  definition was suggested in \cite[5.13]{W-alg}. For monoidal
categories similar conditions were considered by Takeuchi \cite[Definition 5.1]{Tak} and in \cite{M}.
Notice that the term {\em bimonad} is used with a different meaning in
by Brugui\`eres and Virelizier (see \ref{opmon}).

\begin{definition}\label{D.1.5} \em
A {\em bimonad} $\bH$ on a category $\A$ is an endofunctor $H : \A \to \A$
 which has a monad structure
 $\uH=(H, m, e)$ and a comonad structure $\oH=(H, \delta, \varepsilon)$ such that
 \begin{rlist}
\item $\varepsilon : H \to 1$ is a morphism from the  monad $\uH$ to the identity monad;
 \item $e : 1 \to H$ is a morphism from the identity comonad to the  comonad $\oH$;
\item
  there is a mixed distributive law
 $\lambda:HH\to HH$ from the monad $\uH$ to the comonad $\oH$ yielding the commutative
 diagram

\begin{equation}\label{D.1.18}
\xymatrix{
HH \ar[r]^-{m} \ar[d]_-{H \delta}& H \ar[r]^-{\delta}& HH \\
HHH \ar[rr]_-{\lambda H}&& HHH , \ar[u]_-{Hm}} \qquad
\end{equation}
Note that the conditions (i), (ii) just mean commutativity of the
diagrams
\begin{equation}\label{D.1.18b}
 \xymatrix{ HH \ar[r]^{H\varepsilon} \ar[d]_m & H \ar[d]^\varepsilon \\
    H \ar[r]^\varepsilon & 1 , } \qquad
\xymatrix{
1 \ar[r]^-{e} \ar[d]_-{e} & H \ar[d]^-{\delta}\\
H \ar[r]_-{eH}& HH}, \qquad
\xymatrix{ 1\ar[r]^e \ar[dr]_= & H \ar[d]^\vareps\\
         & 1 .}
\end{equation}
\end{rlist}
\end{definition}

\begin{thm}\label{Hopf-mod}{\bf Hopf modules.} \em
Given a bimonad $\bH =(\uH,\oH,\lambda)$ on $\A$,
the objects of $\A^{\oH}_{\uH}(\lambda)$
are called {\em mixed $H$-bimodules} or {\em $H$-Hopf modules}.
By \ref{mon-com}, they are
 triples $(a, h_a, \theta_a)$, where
$(a, h_a) \in \A_{\uH}$ and $(a,\theta_a) \in \A^{\oH}$
with commuting diagram

\begin{equation}\label{H.1}
\xymatrix{
\uH(a) \ar[r]^-{h_a} \ar[d]_-{\uH(\theta_a)}& a \ar[r]^-{\theta_a}& \oH(a) \\
\uH\oH(a) \ar[rr]_-{\lambda_a}&& \oH\uH(a). \ar[u]_-{\oH(h_a)}}
\end{equation}
The morphisms in $\A^{\oH}_{\uH}(\lambda)$ are morphisms in $\A$
which are $\uH$-monad as well as $\oH$-comonad morphisms,
\end{thm}

Recall that a morphism $q: a \to a$ in a category $\A$ is an {\em idempotent} when $qq=q$, and an idempotent $q$ is said to {\em split} if $q$ has a factorization
$q=i \cdot \bar{q}$ with $\bar{q} \cdot i=1$. This happens if and only if
the equaliser $i={\rm Eq}(1_a,q)$ exists or - equivalently - the
coequaliser $\bar q= {\rm Coeq}(1_a,q)$ exists (e.g. \cite[Proposition 1]{BoDe}).
The catgeory $\A$ is called {\em Cauchy complete} provided every idempotent
in $\A$ splits.

\begin{thm}\label{comp-fun}{\bf Comparison functors.} \em
Given a bimonad $\textbf{H}=(\uH=(H, m,e),\oH=(H, \delta, \varepsilon), \lambda)$
on a category $\A$, the mixed distributive law $\lambda$ induces
functors

$$\begin{array}{rl}
K_{\uH} : \A \to (\A_{\uH})^{\woH},&
 a \mapsto ((H(a), m_a), \delta_a), \\[+1mm]
 K_{\oH} : \A \to (\A^{\oH})_{\wuH},&
 a \mapsto ((H(a), \delta_a), m_a),
\end{array}$$
where $\widehat{\oH}$ is the lifting of the comonad $\oH$ and
$\widehat{\uH}$ is the lifting of the  monad $\uH$
by the mixed distributive law $\lambda$. We know that
$(\A_{\uH})^{\woH}\simeq (\A^{\oH})_{\wuH}$ and denote this category by $\A_{\uH}^{\oH}(\lambda)$
(see \ref{mon-com}). There are commutative diagrams

\begin{equation}\label{D.1.20}
\xymatrix{
\A \ar[rr]^-{K_{\uH}} \ar[rrdd]_-{\phi_{\uH}} && (\A_{\uH})^{\widehat{\oH}} \ar[dd]^-{U^{\widehat{\oH}}}\\\\
&& \A_{\uH}, }\quad
\xymatrix{ \A \ar[rr]^{K_\oH} \ar[rrdd]_{\phi^\oH}
&& (\A^\oH)_{\widehat\uH}
\ar[dd]^{U_{\widehat\uH}}\\\\
&& \A^\oH. }
\end{equation}

{\bf (i) The functor $\phi_\uH$.} The forgetful functor $U_{\uH}:
\A_{\uH} \to\A$ is right adjoint to the free functor $\phi_{\uH}$
and the unit $\eta_{\uH}: 1 \to U_{\uH} \phi_{\uH} $ of this
adjunction is the natural transformation $e: 1 \to H$. Since
$\varepsilon: H \to 1$ is a morphism from the monad ${\uH}$ to the
identity monad, $\varepsilon \cdot e=1$, thus $e$ is a split
monomorphism.

The adjunction $\phi_{\uH} \dashv U_{\uH}$ generates
the comonad $\phi_{\uH} U_{\uH}$ on $\A_{\uH}$.
 Recall that for any $(a, h_a) \in \A_{\uH}$,
$\phi_{\uH} U_{\uH}(a, h_a)=(H(a), m_a)$ and
 $\widehat{\oH}(a, h_a)=(H(a), H(h_a) \cdot \lambda_a ).$

 As pointed out in \cite{M}, for any object $b$ of $\A$,
 $K_\uH (b)= (H(b), \alpha_{H(b)})$ for some $\alpha: H(b)\to HH(b)$,
 thus inducing a natural transformation
  $$\alpha_{K_{\uH}}: \phi_{\uH} \to \widehat{\oH} \phi_{\uH},$$
whose component at $b \in \A$ is $\alpha_{H(b)}$, we may choose it to be just
  $\delta_b$,
and we have a morphism of comonads
$$ t_{K_{\uH}}:
\xymatrix{\phi_{\uH} U_{\uH} \ar[rr]^-{\alpha_{K_{\uH}}U_{\uH}} && \widehat{\oH}\phi_{\uH}  U_{\uH} \ar[rr]^{\widehat{\oH}\varepsilon _{\uH} } && \widehat{\oH} ,}$$
where $\varepsilon_{\uH}$ is the counit of the adjunction
$\phi_{\uH}  \dashv U_{\uH}$, and since
$(\varepsilon_{\uH})_{(a,h_a)}=h_a,$ we see that for all $(a, h_a)
\in \A_{\uH}$, $(t_{K_{\uH}})_{(a,h_a)}$ is the composite
\begin{equation}
 \xymatrix{H(a) \ar[r]^-{\delta_a}&HH(a) \ar[r]^-{H(h_a)}& H(a).}
\end{equation}

{\bf (ii) The functor $\phi^\oH$.}
 The cofree $\oH$-comodule functor $\phi^\oH$  has
 the forgetful functor
$U^\oH: \A^\oH \to \A$ as a left adjoint. The
unit $\eta : 1 \to \phi^\oH U^\oH$ and counit
$\sigma: U^\oH\phi^\oH\to 1$ of the adjunction
$U^\oH \dashv \phi^\oH$ are given by the
formulas:
$$\eta_{(a,\, \theta_a)}=\theta_a : (a, \theta_a) \to \phi^\oH
U^\oH(a, \theta_a)=(H(a), \delta_a)$$ and
$$\sigma_a=\varepsilon_a : H(a)=U^\oH\phi^\oH(a) \to a.$$
Since $\varepsilon$ is a split epimorphism, it follows from
Corollary 3.17 of \cite{Me} that, when $\A$ is Cauchy complete, the
functor $\phi^\oH$ is monadic.

Since $K_{\overline{H}}(a)=((H(a), \delta_a), m_a)$, it is easy to
see that the $a$-component of
$$\alpha_{K_{\overline{H}}} :
\widehat{\underline{H}}K_{\overline{H}} \to K_{\overline{H}}$$ is
just the morphism $m_a : HH(a) \to H(a),$ and we have a monad
morphism
$$\xymatrix{
t_{K_{\oH}}: \wuH \ar[rr]^{\widehat{\uH} \eta} & &
 \uH \phi^\oH U^\oH \ar[rr]^{\alpha_{K_{\oH}}U^{\uH}} & & \phi^\oH U^\oH.
}$$
It follows that for any
$(a,\, \theta_a) \in \A^{\overline{H}}\,$,
$(t_{K_{\overline{H}}})_{(a,\, \theta_a)}$  is the composite
\begin{equation}\xymatrix{H(a) \ar[r]^-{H(\theta_a)}& HH(a)
\ar[r]^-{m_a}&H(a)\,.}\end{equation}
\end{thm}

\begin{thm}\label{T.1.8} {\bf The comparison functor as an equivalence.}
 Let $\A$ be a Cauchy complete category. For a bimonad
 $\bH=(\uH=(H, m,e),\oH=(H, \delta, \varepsilon), \lambda)$,
  the following are equivalent:
\begin{blist}
\item
 $K_H: \A \to \A^\oH_\uH(\lambda),\; a \to (H(a), \delta_a, m_a),$
    is an equivalence of categories;
\item $t_{K_{\uH}}:\phi_\uH U_\uH\to\widehat{\oH}$ is an isomorphism of comonads;

\item for any $(a, h_a)\in \A^{\uH}$,
    the composite $H(h_a)\cdot \delta_a$ is an isomorphism;

\item $t_{K_{\oH}}:\widehat{\uH}\to \phi^{\oH}U^{\oH}$ is an isomorphism of monads;

\item  for any $(a,\,\theta_a) \in \A^{\overline{H}}\,$,
 the composite $m_a \cdot H(\theta_a)$ is an isomorphism.
\end{blist}
\end{thm}

\begin{proof} We may identify the functors $K_\uH$, $K_\oH$ and $K_H$.

(a)$\LRa$(b) Since $\A$ is Cauchy complete and since the unit
$\eta_{\uH}: 1 \to U_{\uH} \phi_{\uH} $ of the adjunction
$\phi_{\uH} \dashv U_{\uH}$ is a split monomorphism, the functor
$\phi_{\uH} $ is comonadic by the dual of \cite[Theorem 6]{Mes}.
  Now, by \cite[Theorem 4.4.]{M}, $K_{\uH}$ is an equivalence if and only if
  $t_{K_{\uH}}$ is an isomorphism.
\smallskip

(b)$\LRa$(c) and (d)$\LRa$(e). By \ref{comp-fun}, the morphisms in (b) come out as the morphisms in (c), and the morphisms in (d) are just those in (e).\smallskip
\smallskip

(a)$\LRa$(d) Since $\varepsilon$ is a split epimorphism, it follows from
 \cite[Corollary 3.17]{Me} that (since $\A$ is Cauchy complete) the
functor $\phi^{\overline{H}}$ is monadic and hence $K$ is an equivalence
by \ref{MoB.0.1}.
\end{proof}

\section{Antipode}

We consider
 a bimonad $\bH =(H,m,e,\delta,\varepsilon,\lambda)$ on any catgeory $\A$.

\begin{thm}\label{canon}{\bf Canonical maps.} \em
Define the composites
 \begin{equation}\label{gamma}
\begin{array}{l}
\gamma:\xymatrix{HH \ar[r]^-{\delta H}& HHH \ar[r]^-{Hm}& HH,} \\
   \gamma':\xymatrix{HH \ar[r]^-{H\delta}& HHH \ar[r]^-{mH}& HH.}
\end{array}
\end{equation}
 In  the diagram
 $$\xymatrix{ HHH \ar[rr]^-{\delta HH} \ar[dd]_-{Hm}&& HHHH
\ar[rr]^-{HmH} \ar[dd]^-{HHm}&& HHH
  \ar[dd]^-{Hm}\\\\ HH \ar[rr]_-{\delta H}&& HHH \ar[rr]_-{Hm} &&HH, }$$
the left square commutes by naturality of $\delta$, while the right square commutes by associativity of $m$. From this we see that  $\gamma$ is left $\uH$-linear as a morphism from $(HH, Hm)$ to itself.
A similar diagram shows that $\gamma'$ is right $\uH$-linear
as a morphism from $(HH, mH)$ to itself.
Moreover, in the diagram
$$\xymatrix{ H \ar[rr]^-{He} \ar[dd]_-{\delta}&& HH \ar[rr]^-{\delta H}&& HHH \ar[dd]^-{Hm}\\\\
HH \ar[rrrruu]^-{HHe} \ar@{=}[rrrr]&&&&HH}$$
the top triangle commutes by functoriality of composition, while the bottom triangle commutes because $m \cdot He=1$.
Drawing a similar diagram for $H\delta$ and $mH$, we obtain
\begin{equation}\label{E.1.24}
\gamma \cdot He= \delta, \quad \gamma' \cdot eH= \delta.
\end{equation}
\end{thm}

\begin{definition} \em
A natural transformation $S: H \to H$ is said to be
\begin{bull}
\item a {\em left antipode} if   $m \cdot (SH) \cdot \delta = e \cdot \varepsilon$;
\item a {\em right antipode} if   $m \cdot (HS) \cdot \delta = e \cdot  \varepsilon$;
\item an {\em antipode} if it is a left and a right antipode.
\end{bull}
A bimonad $\bH$ is said to be a {\em Hopf monad} provided it has an antipode.
\end{definition}

Following the pattern of the proof of \cite[15.2]{BW} we obtain:

\begin{proposition}\label{P.1.11}
We refer to the notation in \ref{canon}.
\begin{zlist}
\item
If $\gamma$ has an $\oH$-linear left inverse, then $\bH$ has a left antipode.
\item
If $\gamma'$ has an $\oH$-linear left inverse, then $\bH$ has a right antipode.
\end{zlist}
\end{proposition}

\begin{proof}
(1) Suppose there exists an $\bH$-linear morphism $\beta: HH \to HH$ with $\beta \cdot \gamma=1$.
 Consider the composite
 $$\xymatrix{ S : H \ar[r]^{He} & HH \ar[r]^{\beta} & HH \ar[r]^{\varepsilon H}& H.}$$
 We claim that $S$ is a left antipode of $\bH$. Indeed, in the  diagram
$$\xymatrix{
H  \ar[r]^{\delta}& HH \ar[rr]^{HeH} \ar@{=}[drr]&& HHH \ar@{}[rrd]^{(1)} \ar[rr]^{\beta H}\ar[d]_{Hm}&& HHH \ar@{}[rrd]^{(2)} \ar[d]^{Hm}\ar[rr]^{\varepsilon HH}&&HH \ar[d]^{m}\\
&&& HH \ar[rr]_{\beta}&& HH \ar[rr]_{\varepsilon H}&&H
\,,}$$
 the {\em triangle} commutes since $e$ is the unit for the monad $\uH$,
{\em rectangle} (1) commutes by $\oH$-linearity of $\beta$, and
{\em rectangle} (2) commutes by naturality of $\varepsilon$.
Thus
$$m \cdot SH \cdot \delta=m \cdot \varepsilon HH \cdot \beta H \cdot HeH \cdot \delta=\varepsilon H \cdot\beta \cdot\delta,$$
and using (\ref{E.1.24}), we have
$$\varepsilon H \cdot \beta \cdot \delta=
  \varepsilon H \cdot \beta \cdot \gamma \cdot He=\varepsilon H \cdot He=
   e \cdot \varepsilon.$$
Therefore $S$ is a left antipode of $\bH$.
\smallskip

(2) Denoting the left inverse of $\gamma'$ by $\beta'$,
it is shown along the same lines that $S'=H\ve\cdot \beta'\cdot eH$
is a right antipode.
\end{proof}

\begin{lemma}\label{L.1.12}
Suppose that $\gamma$ is an epimorphism. If $f, g : H \to H $ are two natural transformations such that 
$$m \cdot fH \cdot \delta = m \cdot gH \cdot \delta  \quad \mbox{ or } \quad 
  m \cdot Hf \cdot \delta = m \cdot Hg \cdot \delta,$$ 
then $f=g$.
\end{lemma}

\begin{proof} Assume $m \cdot fH \cdot \delta = m \cdot gH \cdot \delta$.
Since $\gamma \cdot He =\delta $ by (\ref{E.1.24}), we have
$$m \cdot fH \cdot \gamma \cdot He =m \cdot gH \cdot \gamma \cdot He,$$
and, since $\gamma$ is also $\uH$-linear, it follows by
Lemma \ref{L.1.9} that $$m \cdot fH \cdot \gamma = m \cdot gH \cdot \gamma. $$
But $\gamma$ is an epimorphism by our assumption, thus
$$m \cdot fH=m \cdot gH.$$
 By naturality of $e: 1 \to H$, we have the commutative diagrams
$$\xymatrix{
H \ar[d]_{He}\ar[r]^{f}& H \ar[d]^{He}\\
HH \ar[r]^{fH}&HH,} \qquad
\xymatrix{
H \ar[d]_{He}\ar[r]^{g}  & H \ar[d]^{He}\\
HH \ar[r]^{gH}  &HH.
}
$$
Thus, since $m \cdot He=1$,
$$f=m\cdot He\cdot f= m\cdot fH\cdot He=m\cdot gH\cdot He=m\cdot He\cdot g=g.$$

If  $m \cdot Hf \cdot \delta = m \cdot Hg \cdot \delta$ similar arguments apply.
\end{proof}

\begin{thm}\label{P.1.13}{\bf Characterising Hopf monads.} 
Let $\bH =(H,m,e,\delta,\varepsilon,\lambda)$ be a bimonad. 

\begin{zlist}
\item The following are equivalent:
\begin{blist}
\item    $\gamma=Hm \cdot \delta H   : HH \to HH$ is an isomorphism;
\item $\gamma'=mH\cdot H\delta : HH \to HH$ is an isomorphism;
\item $\bH$ has an antipode.
\end{blist}
\item
If $\bH$ has an antipode and $\A$ admits equalisers,
then the comparison functor (see \ref{comp-fun}) 
$$K_{\underline{H}}: \A \to
{\A_{\underline{H}}^{\overline{H}}}(\lambda)$$ makes $\A$
(isomorphic to) a coreflective subcategory of the category
${\A_{\underline{H}}^{\overline{H}}}(\lambda)$.
\end{zlist}
\end{thm}
\begin{proof} (1)
(c)$\Rightarrow$(a) The proof for \cite[Proposition 6.10]{M} applies almost
literally.

(a)$\Rightarrow$(c)
Write $\beta : HH \to HH$ for the inverse of $\gamma$. Since $\gamma$ is $\uH$-linear, it follows that $\beta$ also is $\uH$-linear. Then, by Proposition \ref{P.1.11},
$S= \varepsilon H \cdot \beta \cdot He$ is a left antipode of $\bH$. We show that $S$ is also a right antipode of $\bH$.
In the diagram
$$
\xymatrix{
H \ar[rr]^-{\delta} \ar[ddrr]_-{\delta}&&HH \ar@{}[dd]^{(1)} \ar[rr]^-{\delta H}&&HHH \ar@{}[dd]^{(2)} \ar[rr]^-{HSH}&& HHH \ar@{}[rrdd]^{(3)} \ar[rr]^-{mH} \ar[dd]^-{Hm}&& HH \ar[dd]^-{m}\\\\
&&HH \ar[rruu]^-{H \delta } \ar[rr]^-{H \varepsilon }&& H \ar[rr]^-{He}&& HH \ar[rr]^-{m}&&H\,.
}$$
\begin{bull}
\item (1) commutes by coassociativity of $\delta$;
\item (2) commutes because $S$ is a left antipode of $\bH$;
\item (3) commutes by associativity of $m$.
\end{bull}
Since $m \cdot He=1=m \cdot eH$ and $H \varepsilon \cdot \delta=1=\varepsilon H \cdot \delta,$
it follows that
$$\begin{array}{rl}
m \cdot (m \cdot HS \cdot \delta)H \cdot \delta&=\; m \cdot mH \cdot HSH \cdot \delta H \cdot \delta
 \;=\; m \cdot He \cdot H \varepsilon \cdot \delta\\[+1mm]
 & =\;m \cdot eH \cdot \varepsilon H \cdot \delta\;=\; m \cdot ((e \cdot \varepsilon)H)\cdot \delta.
\end{array}$$
  $\gamma$ being an epimorphism, Lemma \ref{L.1.12} implies
 $m \cdot HS \cdot \delta = e \cdot \varepsilon$,
proving that $S$ is also a right antipode of $\bH$.
 
(b)$\LRa$(c) can be shown in a similar way.
\smallskip

(2) Since
\begin{itemize}
\item to say that $\gamma$ is an isomorphism is to say that
$(t_{K_{\underline{H}}})_{(H(a), m_a)}$ is an isomorphism for all
$a \in \A$;
\item $(H(a), m_a)=\phi_{\underline{H}}(a)$;
\item the unit $\eta_{\underline{H}}: 1 \to \phi_{\underline{H}}
U_{\underline{H}}$ of the adjunction
$\phi_{\underline{H}} \dashv
U_{\underline{H}}$ is just $e: 1 \to H$, which is a split
monomorphism,
\end{itemize} 
we can apply Corollary \ref{C.3.10} to get the desired result.
\end{proof}

Combining    \ref{P.1.13} and  \ref{T.1.8}, we get:

\begin{thm}\label{T.1.14}{\bf Antipode and equivalence.}
Let $\bH=(H, m, e, \delta, \varepsilon, \lambda)$ be a bimonad on a
category $\A$ and assume that $\A$ admits colimits or limits and $H$
preserves them. Then the following are equivalent:
\begin{blist}
\item  $\bH$ has an antipode;
\item  $\gamma=Hm \cdot \delta H   : HH \to HH$ is an isomorphism;
\item  $\gamma'=mH\cdot H\delta : HH \to HH$ is an isomorphism;
\item  $K_H: \A \to \A^\oH_\uH(\lambda),\; a \to (H(a), \delta_a, m_a),$
       is an equivalence.
\end{blist}
\end{thm}

\begin{proof}
 (a)$\LRa$(b)$\LRa$(c) (in any category) is shown in \ref{P.1.13}.

(b)$\LRa$(d)
 Since $H$ preserves
colimits,  the category $\A_{\uH}$  admits colimits and
the functor $U_{\uH} : \A_{\uH} \to \A$ creates them (see, for
example, \cite{Sch}). Thus
\begin{bull}
\item the functor $\phi_{\uH} U_{\uH}$ preserves colimits;
\item any functor $L: \B \to \A_{\uH}$ preserves colimits if
and only if the composite $U_{\uH}L $ does;
so, in particular, the functor $\widehat{\oH}$ preserves colimits,
since $U_{\uH} \widehat{\oH}=H U_{\uH}$ and since the functor $H U_{\uH}$, being the composite of
two colimit-preserving functors, is colimit-preserving.
\end{bull}
The full subcategory of $\A_{\uH}$ given by the free
${\uH}$-modules is dense and since the functors $\phi_{\uH} U_{\uH}$ and
$\widehat{\oH}$ both preserve colimits, it follows from
\cite[Theorem 17.2.7]{Sch} that the natural transformation (see \ref{T.1.8})
$$t_{K_{\uH}}:\phi_{\uH} U_{\uH}\to \widehat{\oH}$$
is an isomorphism if and only if its restriction to the free $\uH$-modules
is so; i.e. if $(t_{K_{\uH}})_{\phi_{\uH}(a)}$ is an isomorphism for
all $a \in \A$. But since $\phi_{\uH}(a)=(H(a), m_a)$, $t_{K_{\uH}}$
is an isomorphism if and only if the composite
$$\xymatrix{ HH(a) \ar[r]^-{\delta_{H(a)}} & HHH(a) \ar[r]^-{H(m_a)}& HH(a)}$$
is an isomorphism for all $a \in \A$, that is, the isomorphism
 $$\gamma:\xymatrix{ HH \ar[r]^-{\delta H} & HHH \ar[r]^-{Hm}& HH}.$$

(c)$\LRa$(d)
 Since the functor $H$ preserves limits, the category
$\A^{\oH}$ admits and the functor $U^{\oH}$
creates limits. Since $\phi^{\oH}$, being right adjoint,
preserves limits, the functor $\phi^{\oH}U^{\oH}$
also preserves limits. Moreover, since the monad
$\widehat\uH$ is a lifting of the monad $\uH$
along the functor $U^{\oH}$,
$U^{\oH}\widehat\uH=\uH U^{\oH}$,
implying that the functor $\widehat\uH$ also preserves
limits. Now, since the full subcategory of $\A^{\oH}$
spanned by cofree $\oH$-comodules is codense, it follows
from the dual of \cite[Theorem 17.2.7]{Sch} that the natural
transformation $t_{K_{\oH}}$ (see \ref{T.1.8}) is an isomorphism if and only if its
restriction to free $\oH$-comodules is so. But for any
$a\in \A$, $(t_{K_{\oH}})_{(H(a), \delta_a)}=m_{H(a)}\cdot
H(\delta_a)$. Thus $t_{K_{\oH}}$ is an
isomorphism if and only if the composite $\gamma'$ is an isomorphism.
\end{proof}

\section{Local prebraidings for Hopf monads}

For any category $\A$
we now fix a system $\bH =(H, m, e, \delta, \varepsilon)$ consisting of an endofunctor
 $H : \A \to \A$ and natural transformations $m:HH \to H$, $e: 1 \to H$, $\delta: H \to HH$ and
 $\varepsilon : H \to 1$ such that the triple $\uH=(H, m, e)$ is a monad and the triple
  $\oH=(H, \delta, \varepsilon)$ is a comonad on $\A$.

\begin{thm}\label{double.ent}{\bf Double entwinings.} \em
A natural transformation $\tau: HH \to HH$ is called a {\em double entwining} if
\begin{rlist}
\item $\tau$ is a mixed distributive law from the monad $\uH$ to the comonad $\oH$;
\item $\tau$ is a mixed distributive law from the comonad $\oH$ to the monad $\uH$.
\end{rlist}
These conditions are obviously equivalent to
\begin{rlist}
\setcounter{rlist}{2}
\item  $\tau$ is a monad distributive law for the monad $\uH$;
\item  $\tau$ is a comonad distributive law for the comonad $\oH$.
\end{rlist}
\end{thm}

  Explicitely (i) encodes the identities
\begin{equation}\label{E.1.4}
He=\tau \cdot eH
\end{equation}
\begin{equation}\label{E.1.5}
H \varepsilon = \varepsilon H \cdot \tau
\end{equation}
\begin{equation}\label{E.1.6}
\delta H \cdot \tau =H \tau \cdot \tau H \cdot H \delta
\end{equation}
\begin{equation}\label{E.1.7}
\tau \cdot mH = Hm \cdot \tau H \cdot H \tau,
\end{equation}
and (ii) is equivalent to the identities
\begin{equation} \label{E.1.8}
eH=\tau \cdot He
\end{equation}
\begin{equation}\label{E.1.9}
\varepsilon H=H \varepsilon  \cdot \tau
\end{equation}
\begin{equation}\label{E.1.11}
H \delta  \cdot \tau = \tau H \cdot H \tau \cdot  \delta H
\end{equation}
\begin{equation}\label{E.1.10}
\tau \cdot Hm = mH \cdot H \tau \cdot  \tau H
\end{equation}

\begin{thm}\label{tau-bim} {\bf $\tau$-bimonad.} \em
Let $\tau: HH \to HH$ be a double entwining. Then $\bH$ is called a
{\em $\tau$-bimonad} provided the
  diagram
\begin{equation}\label{D.1.2}
\xymatrix{
HH \ar[d]_{\delta\delta} \ar[r]^-{m}& H \ar[r]^{\delta}& HH\\
HHHH \ar[rr]_-{H \tau H}&& HHHH \ar[u]_-{mm}}
\end{equation}
is commutative, that is
$$ \delta \cdot  m=mm \cdot H \tau H \cdot \delta \delta=
Hm \cdot mHH \cdot H \tau H \cdot HH \delta \cdot \delta H,
$$
and also the following diagrams commute
\begin{equation}\label{D.1.3}
   \xymatrix{ HH \ar[r]^{H\varepsilon} \ar[d]_m & H \ar[d]^\varepsilon \\
    H \ar[r]^\varepsilon & 1 , }
     \qquad
\xymatrix{
1 \ar[r]^-{e} \ar[d]_-{e}& H \ar[d]^-{\delta}\\
H \ar[r]_-{eH}& HH,}
   \qquad
 \xymatrix{ 1\ar[r]^e \ar[dr]_= & H \ar[d]^\vareps\\
         & 1 .}
\end{equation}
\end{thm}

\begin{proposition}\label{P.1.1}
Let $\bH$ be a $\tau$-bimonad. Then the composite
$$
\xymatrix{\ttau :HH \ar[r]^-{\delta H} & HHH \ar[r]^-{H \tau}& HHH \ar[r]^-{mH}& HH
}$$
is a mixed distributive law from the monad $\uH$ to the comonad $\oH$.
Thus $\bH$ is a bimonad (as in \ref{D.1.5})  with mixed distributive law $\ttau$.
\end{proposition}

\begin{proof}
We have to show that $\ttau$ satisfies
\begin{equation}\label{E.1.12}
He=\ttau \cdot eH
\end{equation}
\begin{equation}\label{E.1.13}
H \varepsilon = \varepsilon H \cdot \ttau
\end{equation}
\begin{equation}\label{E.1.14}
\delta H \cdot \ttau =H \ttau \cdot \ttau H \cdot H \delta
\end{equation}
\begin{equation}\label{E.1.15}
\ttau \cdot mH = Hm \cdot \ttau H \cdot H \ttau
\end{equation}

Consider the   diagram
$$
\xymatrix{
H \ar@{}[rrdd]_-{(1)} \ar[rr]^-{eH} \ar[dd]_-{eH}&& HH \ar@{}[rrdd]_-{(2)} \ar[rr]^-{\tau} \ar[dd]_-{eHH} && HH \ar[dd]_-{eH} \ar@{=}[ddrr] &\\\\
HH \ar[rr]_-{\delta H}&& HHH \ar[rr]^-{H \tau} && HHH \ar[rr]^-{mH} && HH\,, }
$$
which is commutative since
{\em square} (1) commutes by (\ref{D.1.3});
{\em square} (2) commutes by functoriality of composition;
 the {\em triangle} commutes since $e$ is the identity of the monad $\uH$.
Thus $\ttau \cdot eH = mH \cdot H \tau \cdot \delta H \cdot eH = \tau \cdot eH$, and (\ref{E.1.4})
implies $\ttau \cdot eH =He$, showing (\ref{E.1.12}).

Consider now the diagram
$$
\xymatrix{
HH \ar[r]^-{\delta H}&
HHH \ar[r]^-{H \tau} \ar[dr]_-{HH \varepsilon} \ar[dd]_-{\varepsilon HH} &
 HHH \ar@{}[rd]^{(1)} \ar[r]^-{mH} \ar[d]^-{H \varepsilon H}&  HH \ar[d]^-{\varepsilon H}\\
&& HH \ar[r]^-{\varepsilon H}&  H \\
& HH \ar@{}[ru]^-{(2)}\ar[rru]_-{H \varepsilon}  }$$
in which
{\em square} (1) commutes because $\varepsilon$ is a morphism of monads and thus $\varepsilon  \cdot m =\varepsilon \cdot H \varepsilon$;
the {\em triangle} commutes because of (\ref{E.1.5}),
{\em diagram} (2) commutes because of functoriality of composition.

Thus $\varepsilon H \cdot \ttau = \varepsilon H \cdot mH \cdot H \tau \cdot \delta H= H \varepsilon \cdot \varepsilon HH \cdot \delta H = H \varepsilon$, showing (\ref{E.1.13}).

Constructing suitable commutative diagram we can show
$$\begin{array}{rl} 
\ttau \cdot mH & = \;  mH \cdot H \tau \cdot \delta H \cdot mH \\
&= \; mH \cdot HHm \cdot H mHH \cdot HH \tau H \cdot H \tau HH \cdot HHH \tau \cdot \delta \delta H, \\[+2mm]
Hm \cdot \ttau H \cdot H \ttau & =\;
   Hm \cdot m HH \cdot H \tau H \cdot \delta H H
   \cdot HmH \cdot HH \tau \cdot H \delta H \\[+1mm]
&=\; mH \cdot HHm \cdot HmHH \cdot HH \tau H \cdot H \tau HH \cdot HHH\tau \cdot \delta \delta H.
\end{array}
$$
Comparing this two identities we get the condition (\ref{E.1.15}).

To show that (\ref{E.1.14}) also holds, consider the diagram
$$
\xymatrix{
HHH \ar[rr]^-{\delta HH} \ar[rddd]_-{\delta HH} \ar@{}[rrd]^-{(1)} & &
 HHHH    \ar[r]^-{H\tau H} \ar[d]_{H\delta HH}&
HHHH  \ar[rr]^-{mHH} \ar[d]^-{HH\delta H} \ar@{}[rrd]^-{(3)}  &&
  HHH \ar[d]^-{H \delta H}\\
   && HHHHH \ar@{}[r]^-{(2)} \ar[dd]^-{HH\tau H}&  HHHHH
       \ar[dd]^-{HHH \tau} \ar[rr]_-{mHHH}
 \ar@{}[rrdd]^-{(4)}&&  HHHH \ar[dd]^-{HH \tau}\\\\
HH \ar[r]_-{\delta \delta}\ar[uuu]^-{H \delta}& HHHH \ar[ruu]^-{\delta HHH}&HHHHH \ar[ruu]_-{H\tau HH}&  HHHHH \ar[rrd]_-{mmH} \ar[rr]^-{mHHH} && HHHH \ar[d]^-{HmH}\\
&&&&& HHH,}$$
in which
  the {\em triangles} and {\em diagrams} (1) and (3) commute by functoriality
of composition;
{\em diagram} (2) commutes by (\ref{E.1.11});
{\em diagram} (4) commutes by naturality of $m$.

Finally we construct the diagram
$$
\xymatrix{
HH \ar[d]_-{\delta \delta}\ar[rr]^-{\delta H} &&
HHH \ar@{}[rrd]^-{(1)}\ar[lld]^-{HH\delta} \ar[rr]^-{H \tau}&&
HHH \ar@{}[rrddd]^-{(2)} \ar[d]^-{H \delta H} \ar[rr]^-{mH} && HH \ar[d]^-{\delta H}  \\
HHHH \ar@{}[rrdd]^-{(3)} \ar[rr]_{H \tau H} \ar[dd]_-{\delta HHH} && HHHH  \ar@{}[rrdd]^-{(4)} \ar[dd]^-{\delta HH} \ar[rr]^-{HH \tau}&& HHHH \ar[dd]^-{\delta HHH} && HHH \\\\
HHHHH \ar[rr]_-{HH \tau H} && HHHHH \ar[drr]_-{H \tau HH} \ar[rr]^{HHH \tau} &&
HHHHH \ar[rr]^-{H \tau HH} \ar@{}[d]^-{(5)} && HHHHH \ar[uu]^-{mm H}\\
&&&& HHHHH \ar[rru]_-{HHH \tau} }$$
in which
{\em diagram} (1) commutes by (\ref{E.1.6});
{\em diagram} (2) commutes by (\ref{D.1.2}) because $\delta HHH \cdot H \delta H =\delta \delta H$;
 the {\em triangle} and {\em diagrams} (3), (4) and (5) commute by functoriality
of composition.

It now follows from the commutativity of these diagrams that
$$\begin{array}{rl}
\delta H \cdot \ttau& =\; \delta H \cdot mH \cdot H \tau \cdot \delta H \\[+1mm]
 & =\; mmH \cdot HHH \tau \cdot H \tau HH \cdot HH \tau H \cdot \delta HHH \cdot \delta \delta  \\[+1mm]
 & =\; (H mH \cdot HH \tau \cdot H \delta H) \cdot
   (mHH \cdot H \tau H \cdot \delta HH) \cdot H \delta \\[+1mm]
 & =\; H \ttau \cdot \ttau H \cdot H \delta.
\end{array}$$
Therefore $\ttau$ satisfies the conditions (\ref{E.1.12})-(\ref{E.1.15}) and hence is a mixed distributive law from the monad $\uH$ to the comonad $\oH$.
\end{proof}

\begin{corollary}
In the situation of the previous proposition, if
$(a, \theta_a) \in \A^{\oH}$, then
$(H(a), \theta_{H(a)}) \in \A^{\oH}$, where
$\theta_{H(a)}$ is the composite
$$\xymatrix{ H(a) \ar[r]^-{H(\theta_a)}& HH(a) \ar[r]^-{\delta_{H(a)}}& HHH (a)
\ar[r]^-{H\tau_{a}}& HHH(a) \ar[r]^-{m_{H(a)}}& HH (a)}.$$
\end{corollary}

\begin{proof}
Write $\widehat{\uH}$ for the monad on the category $\A^{\oH}$ that is the
lifting of $\uH$ corresponding to the mixed distributive law $\ttau$.
Since $\theta_{H(a)}= \ttau_a \cdot H(\theta_a)$, it follows that $(H(a), \theta_{H(a)})=\widehat{\uH}(a, \theta_a),$ and thus $ (H(a), \theta_{H(a)})$
is an object of the category $\A^{\oH}$.
\end{proof}

\begin{thm}\label{bimod}{\bf $\tau$-Bimodules.}  \em
Given the conditions of Proposition \ref{P.1.1},
we have the commutative diagram (see (\ref{D.1.18}))
$$
\xymatrix{
HH \ar[r]^{m} \ar[d]_-{H \delta}& H \ar[r]^{\delta}& HH \\
HHH \ar[rr]_{\ttau H}&& HHH , \ar[u]_{Hm}}
$$
and thus $H$ is a bimonad by the entwining $\ttau$ and the mixed
bimodules are objects $a$ in $\A$ with a module structure $h_a:H(a)\to a$ and a
comodule structure $\theta_a: a\to H(A)$ with a commutative diagram
$$\xymatrix{
 H(a) \ar[d]_{H(\theta_a)} \ar[r]^{h_a}& a \ar[r]^{\theta_a} & H(a)  \\
HH(a)\ar[rr]^{\ttau_a}   && HH(a) \ar[u]_{H(h_a)}.
}$$
By definition of $\ttau$, commutativity of this diagram is equivalent to the
commtativity of
\begin{equation}\label{D.1.17}
\xymatrix{
& H(a) \ar[dl]_-{H(\theta_a)} \ar[r]^-{h_a}& a \ar[r]^-{\theta_a}&H(a) &\\
HH(a)\ar[rd]_{\delta_{H(a)}} &&&& HH(a) \ar[lu]_{H(h_a)}\\
& HHH(a) \ar[rr]_{H(\tau_a)}&& HHH(a) \ar[ru]_{m_{H(a)}}&.}
\end{equation}  A morphism $f:(a, h_a, \theta_a) \to (a', h_{a'}, \theta_{a'})$ is a morphism $f: a \to a'$ such that $f \in \A^{\oH}$ and $f \in \A_{\uH}$.

We denote the category $\A^\oH_\uH(\ttau)$ by $\A^\nH_\nH$.
\end{thm}
\medskip

\begin{thm}\label{bim.anti}{\bf Antipode of a $\tau$-bimonad.}
Let $\bH=(H,m,e,\delta,\vareps)$ be a $\tau$-bimonad with an antipode $S$
where $\tau:HH\to HH$ is a double entwining. Then
\begin{equation}\label{bim.S.E}
S\z m = m \cdot SS \cdot \tau  \; \mbox{ and }\; \delta\z S= \tau\z SS\z \delta.
\end{equation}

If $\tau \cdot HS=SH \cdot \tau$
and $\tau \cdot SH=HS\cdot \tau$,
then $S:H\to H$ is a monad as well as a comonad morphism.
\end{thm}

\begin{proof}
   Since $(HH, H \tau H \cdot \delta, \varepsilon \varepsilon)$ is a comonad and
 $(H,m,e)$ is a monad, the collection $\Nat(HH,H)$ of all natural transformations
 from $HH$ to $H$ forms a semigroup with unit $e \cdot \varepsilon \varepsilon$
and with product
$$\xymatrix{f*g : HH \ar[r]^{\delta \delta}& HHHH \ar[r]^{H \tau H}& HHHH \ar[r]^-{fg}& HH \ar[r]^{m}&H} .$$
Consider now the diagram
$$ \xymatrix{
&& HH \ar[d]_{m}\ar[dll]_{H \varepsilon}\ar[rr]^{\delta \delta}&& HHHH \ar@{}[dd]^{(2)}\ar[rr]^{H \tau H} && HHHH \ar[d]^{mHH}\\
H \ar@{}[rr]_{(1)}\ar[rrdd]_{\varepsilon}&& H \ar[dd]_{\varepsilon} \ar[drr]_{\delta} &&&&HHH \ar[d]^{SHH}\ar[lld]_{Hm}\\
&&& &HH \ar@{}[d]_{(3)}\ar@{}[rr]^{(4)}\ar[drr]_{SH}&&HHH \ar[d]^{Hm}\\
&& I \ar[rr]_{e}&&H && HH \ar[ll]^{m}}$$
in which the diagrams (1),(2) and (3) commute because $H$ is a bimonad, while diagram (4) commutes by naturality. It follows that
$$m \cdot Hm \cdot SHH \cdot mHH \cdot H \tau H \cdot \delta \delta=e \cdot \varepsilon \cdot H \varepsilon=\varepsilon \varepsilon \cdot e. $$
Thus $S\z m=m^{-1}$ in $\Nat(HH,H)$.
Furthermore, by (a somewhat tedious) computation
 we can show
$$
m\cdot Hm \cdot HHS \cdot HSH \cdot H \tau \cdot mHH \cdot H \tau H \cdot
 \delta\delta =  e\cdot \varepsilon \cdot H \varepsilon=
  e \cdot \varepsilon \varepsilon.
$$

This shows that $m \cdot SS \cdot \tau=m^{-1}$ in $\Nat(HH,H)$.
Thus $m \cdot SS \cdot \tau=S\z m$.
\smallskip

To prove the formula for the coproduct
consider $\text{Nat}(H, HH)$ as a monoid with unit
$ee \cdot \varepsilon$ and the
convolution product
for $f,g \in \text{Nat}(H, HH)$ given by
$$\xymatrix{
  f *g : H \ar[r]^{\delta}&
HH \ar[r]^{fH}& HHH \ar[r]^{HHg}& HHHH \ar[r]^{mm}&HH\,.}$$

By computation we get
$$\begin{array}{rcl}
(\delta \cdot S)*\delta&=&eH \cdot e \cdot \varepsilon = ee \cdot \varepsilon, \\
\delta *(\tau \cdot SS \cdot \delta)
&=& He \cdot e \cdot \varepsilon = ee \cdot \varepsilon.
\end{array}$$
Thus $(\delta \cdot S)*\delta=1$
and $\delta *(\tau \cdot SS \cdot \delta)=1$, and hence $\delta \cdot S=\tau \cdot SS \cdot \delta.$
\smallskip

Now assume
 $\tau \cdot HS=SH \cdot \tau$ and $\tau \cdot SH=HS \cdot \tau.$ Then we have
$$SS \cdot \tau =SH \cdot HS\cdot \tau =SH \cdot \tau \cdot SH=\tau \cdot HS \cdot SH=\tau \cdot SS, \mbox{ thus } $$
$$S\cdot m=m \cdot SS \cdot \tau=m \cdot \tau \cdot SS=m' \cdot SS.$$
Moreover, since $m \cdot He=1,$ we have
$$S \cdot e=m \cdot He \cdot S \cdot e\stackrel{\text{nat}}=m \cdot SH \cdot He \cdot e\stackrel{(\ref{D.1.3})} =m \cdot SH \cdot \delta \cdot e
\stackrel{\text{antip.}}=e \cdot \varepsilon \cdot e
\stackrel{(\ref{D.1.3})}=e\,.$$
Hence $S$ is a monad morphism from $(H,m,e)$ to $(H, m \cdot \tau, e)$.
\smallskip

For the coproduct,  $SS \cdot \tau=\tau \cdot SS$ implies
  $$\delta \cdot S=\tau \cdot SS \cdot \delta=
SS \cdot \tau \cdot \delta=SS \cdot \delta'.$$
Furthermore,
$$\varepsilon \cdot S=\varepsilon \cdot S \cdot H \varepsilon \cdot \delta  \stackrel{\text{nat}}=
\varepsilon \cdot H \varepsilon \cdot SH \cdot \delta
\stackrel{(\ref{D.1.3})}=
\varepsilon \cdot m \cdot SH \cdot \delta
\stackrel{\text{antip.}}=
\varepsilon \cdot e \cdot \varepsilon
\stackrel{\text{(\ref{D.1.3})}}=
\varepsilon.$$
This shows that $S$ is a comonad morphism from  $(H, \delta, \varepsilon)$
to  $(H, \tau \cdot \delta, \varepsilon).$
 \end{proof}
\medskip

It is readily checked that for a bimonad $H$, the composite $HH$ is again
a comonad as well as a monad. However, the compatibility between these
two structures needs an additional property of the double entwining $\tau$.
This will also help to construct a bimonad "opposite" to $H$.

\begin{thm}\label{Yang}{\bf Local prebraiding.} \em
Let  $\tau:HH\to HH$ be a natural transformation.
 $\tau$ is said to satisfy the
{\em Yang-Baxter equation (YB)} if it induces commutativity of the diagram
$$\xymatrix{ HHH \ar[r]^{\tau
H} \ar[d]_{H \tau} & HHH \ar[r]^{H \tau}& HHH \ar[d]^{\tau H}\\
HHH \ar[r]_{\tau H}& HHH \ar[r]_{H \tau}& HHH\, . }
$$

$\tau$ is called a {\em local prebraiding} provided it is a double entwining
(see \ref{double.ent})
and  satisfies the Yang-Baxter equation.
\end{thm}

\begin{thm}\label{HH-bim}{\bf Doubling a bimonad.}
Let $\bH=(H,m,e,\delta,\vareps)$ be a $\tau$-bimonad
where $\tau:HH\to HH$ is a local prebraiding. Then
  $\HH=(HH,\bar m,\bar e,\bar\delta,\bar\vareps) $ is a $\bar\tau$-bimonad with
 $\bar e= ee$, $\bar\vareps=\vareps\vareps$,
\begin{center}
$\bar m:\xymatrix{HHHH\ar[r]^{H \tau H} & HHHH \ar[r]^{\quad mm}&HH\, ,}$ \\[+1mm]
 $\bar\delta:
\xymatrix{ HH \ar[r]^{\delta\delta\quad}& HHHH \ar[r]^{H \tau H}& HHHH }$
\end{center}
and double entwining
$$\bar\tau:\xymatrix{HHHH \ar[r]^{H \tau H}& HHHH
\ar[r]^{\tau HH}& HHHH \ar[r]^{HH \tau}& HHHH\ar[r]^{H \tau H}& HHHH\, .}$$
\end{thm}

\begin{proof}
  We already know that $(HH, \bar{m}, \bar{e})$ is a monad and that
$(HH, \bar{\delta}, \bar{\varepsilon})$ is a comonad. First
we have to show that $\bar{\tau}$ is a mixed distributive law from the monad $(HH, \bar{m}, \bar{e})$ to the comonad$(HH, \bar{\delta}, \bar{\varepsilon})$, that is
$$HH \bar{e}=\bar{\tau} \cdot \bar{e}HH, \quad
HH \bar{\varepsilon}=\bar{\varepsilon}HH \cdot \bar{\tau},$$
$$HH\bar{m}\cdot \bar{\tau}HH \cdot HH\bar{\tau}= \bar{\tau} \cdot \bar{m}HH,$$
$$HH\bar{\tau}\cdot \bar{\tau} HH \cdot HH \bar{\delta}=\bar{\delta}HH\cdot\bar{\tau} .$$

The first two equalities can be verified by placing the composites in
suitable commutative diagrams.
The second two identities are obtained by lengthy standard computations
(as known for classical Hopf algebras).

It remains to show
 that $(HH, \bar{m},\bar{e},\bar{\delta},\bar{\varepsilon})$ satisfies the
 conditions of Definition \ref{D.1.5} with respect to $\bar{\tau}$.
Again

$$\begin{array}{rcl}
\bar{\varepsilon}\cdot \bar{m}&=&
\varepsilon \cdot H \varepsilon  \cdot HH\varepsilon \varepsilon\, =\, \bar{\varepsilon} \cdot HH \bar{\varepsilon}, \mbox{ and }\\[+1mm]
\bar{\delta}\cdot \bar{e}&=&HHee \cdot He \cdot e=HH \bar{e}\bar{e}
\end{array}$$
are shown by standard computations and
$$\bar{\varepsilon}\bar{e}=\varepsilon \cdot \varepsilon H \cdot eH \cdot e \stackrel{\text{(\ref{D.1.18b})}} = \varepsilon \cdot e=1.$$

To show that $(HH, \bar{m},\bar{e},\bar{\delta},\bar{\varepsilon},\bar{\tau})$ satisfies
(\ref{D.1.18}), consider the diagram
\begin{small}
$$\xymatrix{
H^4 \ar@{}[rd]^{(1)}\ar[dd]^{H^2 \delta H} \ar[r]^{H \tau H}
& H^4 \ar[drr]_{H \delta H^2}\ar[rrrr]^{m H^2}
&&&& H^3\ar[drr]_{\delta HH} \ar[rr]^{H m}
&& H^2 \ar@{}[d]^{(3)}\ar[rr]^{\delta H}&&H^3\ar[dd]^{H^2 \delta}\\
&&& H^5 \ar@{}[d]^{(4)}\ar@{}[rr]^{(2)} \ar[drr]^{\tau H^3}&&&& H^4\ar[rru]_{H^2 m} \ar[d]^{H^3 \delta}&& \\
H^5 \ar@{}[rd]_{(5)} \ar[dd]^{H^4 \delta}\ar[r]^{H \tau H^2} &H^5 \ar@{}[rrdd]_{(6)}\ar[rru]^{H^2 \tau H}\ar[dd]_{H^4 \delta}\ar[rr]_{\tau H^3}
&& H^5 \ar@{}[rrdd]_{(7)}\ar[dd]_{H^4 \delta}\ar[rr]_{H^2 \tau H}&&
H^5 \ar@{}[rr]_{(8)} \ar[urr]^{H m H^2}
\ar[dd]_{H^4 \delta}&& H^5 \ar@{}[rr]^{(9)} \ar@{}[dd]_{(10)}
\ar[drr]_{H^2 \tau H^2}
&& H^4 \ar[rr]^{H \tau H}&& H^4\\
&&&&&&&&&H^5 \ar@{}[rru]_{(11)} \ar@{}[rrd]_{(12)}\ar[rr]_{H\tau H^2} \ar[u]_{H^3 m}&&H^5 \ar[u]^{H^3 m}\\
H^6 \ar[dd]^{H^3 \tau H}\ar[r]_{H\tau H^3}
& H^6 \ar@{}[rrdd]^{(14)}\ar[dd]^{H^3 \tau H}\ar@{}[dd]_{(13)} \ar[rr]_{\tau H^4}
&& H^6 \ar@{}[rrrrdd]_{(15)} \ar[dd]^{H^3 \tau H} \ar[rr]_{H^2 \tau H^2}
&& H^6 \ar[uurr]^{HmH^3} \ar[rr]_{H^3 \tau H}
&& H^6 \ar[urr]_{HmH^3} \ar[rr]_{H^2 \tau H^2}
&& H^6 \ar@{}[dd]^{(16)}\ar[rr]_{H \tau H^3}&& H^6 \ar[u]^{H^2 m H^2}\\\\
H^6 \ar[r]_{H \tau H^3} & H^6 \ar[rr]_{\tau H^4}
&& H^6\ar[rrrr]_{H^2 \tau H^2}
&&&&H^6 \ar[rruu]_{H^3 \tau H} \ar[rr]_{H \tau H^3}&&H^6\ar[rruu]_{H^3 \tau H},
}$$
\end{small}
in which
  {\em diagram} (1) commutes because $\tau$ is a mixed distributive law and thus
  $$H\tau \cdot \tau H \cdot H \delta =\delta H \cdot \tau;$$
  the {\em diagrams} (2) and (9) commute by (\ref{D.1.18});
  the {\em diagrams} (3)-(8), (10), (11), (13), (14) and (16) commute by naturality;
  {\em diagram} (12) commutes because $\tau$ is a mixed distributive law
(hence
$H m \cdot \tau H \cdot H \tau=\tau \cdot m H$);
  {\em diagram} (15) commutes by \ref{Yang}.
By commutativity of the whole diagram,
$$\begin{array}{l}
\bar{\delta}\cdot \bar{m}\,=\, H \tau H \cdot H^2 \delta \cdot \delta H \cdot Hm \cdot mH^2 \cdot H \tau H \\[+1mm]
\,=\, H^2 m \cdot H^2 m H^2 \cdot H^3 \tau H \cdot H \tau H^3 \cdot H^2 \tau H^2 \cdot \tau H^4 \cdot H \tau H^3 \cdot H^3 \tau H \cdot H^4 \delta \cdot H^2 \delta H \\[+1mm]
 \,=\, HH \bar{\delta} \cdot \bar{\tau}HH \cdot HH \bar{m},
\end{array}$$
and hence $\HH=(HH,\bar{m},\bar{e},\bar{\delta},\bar{\varepsilon})$ is a $\bar{\tau}$-bimonad.
\end{proof}

\begin{thm}\label{YB1}{\bf Opposite monad and comonad.}
Let $\tau:HH\to HH$ be a natural transformation
satisfying the Yang-Baxter equation.
\begin{zlist}
\item If $(H,m,e)$ is a monad and $\tau$ is monad distributive,
then $(H,m\cdot \tau,e)$ is also a monad and $\tau$ is monad distributive for it.

\item If $(H,\delta,\vareps)$ is a comonad and $\tau$ is comonad distributive,
then $(H,\tau\cdot \delta,\vareps)$ is also a comonad and $\tau$ is comonad
distributive for it.
\end{zlist}
\end{thm}

\begin{proof}
(1) To show that $m\cdot \tau$ is associative
construct the diagram
$$\xymatrix{HHH \ar[rr]^{\tau_H}\ar[dd]_{H_\tau}\ar@{}[ddrr]^{(1)} & &
    HHH \ar[r]^{m_H} \ar[d]_{H\tau} \ar@{}[ddr]^{(2)}
& HH \ar[dd]^\tau \\
&& HHH  \ar[d]_{\tau_H}    \\
HHH \ar[r]^{\tau_H} \ar[d]_{Hm} \ar@{}[drr]^{(3)} & HHH \ar[r]^{H_\tau}  &
HHH \ar[r]^{Hm} \ar[d]_{mH} \ar@{}[dr]^{(4)}   & HH \ar[d]^m \\
HH \ar[rr]^\tau    & &HH \ar[r]^m & H ,}$$
where the {\em rectangle} (1) is commutative by the YB-condition,
(2) and (3) are commutative by the monad distributivity of $\tau$, and
the {\em square} (4) is commutative by associativity of $m$.
Now commutativity of the outer diagram shows associativity of
$m\cdot \tau$.

From \ref{lift.mon.mon} we know that $\tau\cdot e_H = He$ and $\tau\cdot He = e_H$
and this implies that $e$ is also the unit for $(H,m\cdot \tau,e)$.

The two pentagons for monad distributivity of $\tau$ for $(H,m\cdot m,e)$ can be read
from the above diagram
by combining the two top rectangles as well as the two left hand rectangles.

\smallskip

(2) The proof is dual to the proof of (1).
\end{proof}

\begin{thm}\label{YB2}{\bf Opposite bimonad.}
Let $\tau:HH\to HH$ be a local prebraiding with $\tau^2=\id$ and
let $\bH=(H,m,e,\delta,\vareps)$ be a $\tau$-bimonad on $\A$.
 Then:
\begin{zlist}
\item $\bH'=(H, m\cdot \tau, e, \tau\cdot \delta, \vareps)$ is also
      a $\tau$-bimonad.
\item If $\bH$ has an antipode $S$ with $\tau \cdot HS=SH \cdot \tau$
and   $\tau \cdot SH=HS\cdot \tau$,
then  $S$ is a $\tau$-bimonad morphism
between the $\tau$-bimonads $\bH$ and $\bH'$.

In this case $S$ is an antipode for $\bH'$.
\end{zlist}
\end{thm}

\begin{proof}
(1) By (1), (2) in \ref{YB1}, $\tau$ is a (co)monad distributive law
from the (co)monad $H$ to the (co)monad $H'$,
and $\varepsilon' \cdot e'=\varepsilon \cdot e=1$ by (\ref{D.1.3}).
Moreover,
$$\varepsilon' \cdot m'=\varepsilon \cdot m \cdot \tau
\stackrel{(\ref{D.1.3})}{=}\varepsilon \cdot H \varepsilon \cdot \tau
 \stackrel{\ref{common}}=\varepsilon \cdot \varepsilon H=\varepsilon \cdot H \varepsilon=\varepsilon' \cdot H \varepsilon',\;\mbox{ and }$$
$$\delta' \cdot e'=\tau \cdot \delta \cdot e
\stackrel{(\ref{D.1.3})}=\tau \cdot eH \cdot e
\stackrel{\ref{mon-com}}=He \cdot e=eH \cdot e=e' H \cdot e'.$$
To prove compatibility for $\bH '$ we have to show the commutativity of the diagram
\begin{equation}\label{opp}
\xymatrix{
HH \ar[r]^{m'} \ar[d]_{\delta' \delta'} & H \ar[r]^{\delta'}&HH\\
HHHH \ar[rr]_{H \tau H}&& HHHH \ar[u]_{m'm'}.}
\end{equation}
For this standard computations (from Hopf algebras) apply.

(2) By \ref{bim.anti},
  $S$ is a $\tau$-bimonad morphism from the $\tau$-bimonad $\bH $ to the
 $\tau$-bimonad $\bH '$.
To show that $S$ is an antipode for $\bH '$ we need the equalities
\begin{center}
 $m' \cdot SH \cdot \delta'=e' \cdot \varepsilon'=e \cdot \varepsilon$
\quad and \quad $m' \cdot HS \cdot \delta'=e' \cdot \varepsilon'=e \cdot \varepsilon.$
\end{center}
Since $\tau \cdot SH=HS \cdot \tau$, we have
$$m' \cdot SH \cdot \delta'=m \cdot \tau \cdot SH \cdot \tau \cdot \delta
=m \cdot HS \cdot \tau \cdot \tau \cdot \delta \stackrel{\tau^2=1}{=}m \cdot HS  \cdot \delta = e \cdot \varepsilon.$$
Since $\tau \cdot HS=SH \cdot \tau$, we have
$$m' \cdot HS \cdot \delta'=m \cdot \tau \cdot HS \cdot \tau \cdot \delta=m \cdot SH \cdot \tau \cdot \tau \cdot \delta\stackrel{\tau ^2=1}=m \cdot SH  \cdot \delta=e \cdot \varepsilon.$$\end{proof}

As we have seen in Theorem \ref{T.1.14}, the existence of an antipode
for a bimonad $\bH$ on a category $\A$ is
equivalent to the comparison functor being an equivalence
provided $\A$ is Cauchy complete and $H$ preserves
colimits. Given a local prebraiding the latter condition on $H$
can be replaced by conditions on the antipode
(compare \cite[Theorem 3.4]{BeDr}, \cite[Lemma 4.2]{BeKe}
 for the situation in braided monoidal category).

\begin{thm}\label{id.split} {\bf Antipode and equivalence.}
 Let $\tau :HH \to HH$ be a local prebraiding and let
$\bH =(H, m, e, \delta, \varepsilon)$ be a $\tau$-bimonad on a
category $\A$ in which idempotents split. Consider the category of
bimodules
 $$ \A^{\nH}_{\nH}= \A^{\oH}_{\uH}(\ttau),$$
where $\ttau=mH \z H \tau \z \delta H $ (see \ref{P.1.1}, \ref{bimod}).

If $\bH $ has an antipode $S$ such that
$\tau \z SH=HS \z \tau$ and $\tau \z HS=SH \z \tau$, then
 the comparison functor
 $K_{\uH}: \A \to \A^{\nH}_{\nH}$
is an equivalence of categories.
\end{thm}

\begin{proof}
 We know that the functor $K_{\uH}$ has a right adjoint if for each
$(a, h_a, \theta_a)\in \A^{\nH}_{\nH},$
the equaliser of the $(a, h_a, \theta_a)-$component of the pair of functors
\begin{equation}\label{(1)}
\xymatrix{ U_{{\uH}} U^{\widehat{\oH}}
\ar@{->}@<0.5ex>[rr]^-{U_{{\uH}} U^{\widehat{\oH}} e_{\widehat{\oH}}} \ar@ {->}@<-0.5ex>
[rr]_-{\beta_{U_{\uH}} U^{\widehat{\oH}}}&& U_{{\uH}} \widehat{\oH} U^{\widehat{\oH}}=U_{{\uH}} U^{\widehat{\oH}} \phi^{\widehat{\oH}} U^{\widehat{\oH}}}
\end{equation}
exists. Here  $e_{\widehat{\oH}}: 1 \to \phi^{\widehat{\oH}} U^{\widehat{\oH}}$ is the unit of the adjunction $U^{\widehat{\oH}} \dashv \phi ^{\widehat{\oH}}$ and $\beta_{U_{\uH}} $ is the composite $$\xymatrix{U_{{\uH}} \ar[rr]^-{e_{{\uH}} U_{{\uH}}}&&U_{{\uH}} \phi_{\uH}U_{{\uH}} \ar[rr]^{U_{\uH}(t_{K_{\uH}})}&& U_{\uH} \widehat{\oH}. }$$
Using the fact that for any $(a, h_a)\in \A_{\uH} $, $$(t_{K_{\uH}})_{(a, h_a)}=H(h_a) \z \delta_a \quad \mbox{ and }$$
$$H(h_a) \z \delta_a \z e_a=H(h_a) \z H(e_a) \z e_a=e_a, $$
it is not hard to show that the $(a, h_a, \theta_a)$-component of
Diagram \ref{(1)} is the pair
$$\xymatrix{a \ar@ {->}@<-0.5ex>[r]_-{\theta_a} \ar@{->}@<0.5ex>[r]^-{e_a}& H(a).}$$
Thus, $K_{\uH}$ has a right adjoint if for each $(a, H_a, \theta_a)\in
\A^{\nH}_{\nH},$ the equaliser of the pair of morphisms $(e_a, \theta_a)$ exists.

Suppose now that $\bH $ has an antipode $S : H \to H$. For each $(a,
h_a, \theta_a)\in \A^{\nH}_{\nH},$ consider the composite $q_a=h_a
\z S_a \z \theta_a: a \to a$. By a (tedious) standard computation -
applying \ref{D.1.17}, \ref{bim.anti}, \ref{common} - one can show
\begin{center}
 $e_a \z q_a= \theta_a \z q_a$ and $q_a \z q_a=q_a $.
\end{center}

\begin{remark}\label{rem.proof} \em Dually, one can prove that for each $(a, h_a, \theta_a)\in \A^{\nH}_{\nH},$ $q_a \z \varepsilon_a =q_a \z h_a,$ thus
$i_a \z \bar q_a \z \varepsilon_a =i_a \z \bar q_a \z h_a$, and since $i_a$ is a (split) monomorphism, it follows that\;
$\q_a \z \varepsilon_a = \q_a \z h_a$.
\end{remark}
Since idempotents split in $\A$, there exist morphisms  $i_a : \bar{a} \to a$ and $\q_a : a \to \bar{a}$ such that $\q_a \z i_a=1_a$ and $i_a \z \q_a=q_a.$ Since $\q_a$ is a (split) epimorphism and since $e_a \z i_a \z \q_a=e_a \z q_a=\theta_a \z q_a=\theta \z i_a \z \q_a,$ it follows that
\begin{equation}\label{(3)} e_a \z i_a=\theta_a \z i_a . \end{equation}
Using this equality it is straightforward to show that the
 diagram
\begin{equation}\label{(4)}
\xymatrix{
\bar{a} \ar[rr]_{i_a} && a \ar@/_/[ll]_{\q_a} \ar@{->}@<0.9ex>[rr]_-{e_a} \ar@ {->}@<-0.9ex>
[rr]_-{\theta_a} && H(a) \ar@/_/@<-1.5ex>[ll]_{h_a \z S_a}}\end{equation} is a split equaliser.
%
Hence for any $(a, h_a, \theta_a)\in \A^{\nH}_{\nH},$ the equaliser
of the pair of morphisms $(e_a, \theta_a)$ exists, which implies
that the functor $K_{\uH}$ has a right adjoint $R_{\uH}:
\A^{\nH}_{\nH} \to \A $ which is given by $$R_{\uH}(a, H_a,
\theta_a)=\bar{a}.$$

Since for any $(a, h_a, \theta_a)\in \A^{\nH}_{\nH},
$\begin{itemize}
  \item $\delta_a \z e_a=e_{H(a)} \z e_a$ \mbox{ and }
        $\varepsilon_a \z e_a=1$ by \ref{tau-bim};
  \item $\varepsilon_{H(a)} \z \delta_a=1$, since $(H, \varepsilon, \delta)$
        is a comonad;
  \item $e_a \z \varepsilon_a =\varepsilon_{H(a)} \z e_{H(a)}$ by naturality,
\end{itemize}
we get a split equaliser  diagram
$$\xymatrix{
a \ar[rr]_{e_a} && H(a) \ar@/_/[ll]_{\varepsilon_a} \ar@{->}@<0.9ex>[rr]_-{e_{H(a)}} \ar@ {->}@<-0.9ex>
[rr]_-{\delta_a} && H^2(a) \ar@/_/@<-1.5ex>[ll]_{H(\varepsilon_a)}}.$$
 This is preserved by any functor, and since $R_{\oH}(H(a), m_a, \delta_a)$
 is the equaliser of the pair of morphisms $(e_{H(a)}, \delta_a)$,
 in particular  $a \simeq R_{\oH}(H(a), m_a, \delta_a)= R_{\oH}(K_{\oH}(a))$.
 Thus $R_{\oH}K_{\oH}\simeq 1$.

For any $(a, h_a,\theta_a)\in\A^{\nH}_{\nH},$ write $\alpha_a$ for the composite $h_a \z H(i_a):H(\bar{a}) \to a.$ We claim that $\alpha_a$ is a morphism in $\A^{\nH}_{\nH}$ from $K_{\uH}(\bar{a})=(H(\bar{a}), m_{\bar{a}}, \delta_{\bar{a}} )$ to $(a, h_a, \theta_a)$.
Indeed, we have
$$\begin{array}{rcl}
\alpha_a \z m_{\bar{a}}&=&h_a \z H(i_a) \z m_{\bar{a}} \\[+1mm]
 \text{\scriptsize{naturality}} &=& h_a \z m_a \z H^2(i_a)\\[+1mm]
\text{\scriptsize{
 $ (a, h_a)\in \A_{\uH}$}} &=&
  h_a \z H(h_a) \z H^2(i_a)=h_a \z H(H(h_a) \z i_a)=h_a \z H(\alpha_a),
\end{array}$$
and this just means that $\alpha_a$ is a morphism in $\A_{\uH}$ from $(H(\bar{a}), m_{\bar{a}})$ to $(a, h_a)$.

Next - using \ref{D.1.17}, \ref{(3)} - we compute
$$ \theta_a \z \alpha_a=\,H(\alpha_a) \z \delta_{\bar{a}}.$$

Thus, $\alpha_a$ is a morphism in $\A^{\oH}$ from $(H(\bar{a}), \delta_{\bar{a}})$ to $(a, \delta_a )$, and hence $\alpha_a$ is a morphism in $\A^{\nH}_{\nH}$ from $K_{\uH}(\bar{a})=(\bar{a}, m_{\bar{a}}, \delta_{\bar{a}} )$ to $(a, h_a, \theta_a)$.

Similarly it is proved that the composite $\beta_a=H(\q_a)\z \theta_a: a \to H(\bar{a})$ is a morphism in $\A^{\oH}$ from $(a, h_a, \delta_a )$ to
$(H(\bar{a}), m_{\bar{a}}, \delta_{\bar{a}})$
and a further calculation yields
\begin{center}
 $\alpha_a \z \beta_a=1_a$ and $\beta_a \z \alpha_a=1_{H(\bar{a})}$.
\end{center}

Hence we have proved that for any $(a, h_a, \theta_a) \in \A^{\nH}_{\nH},$ $\alpha_a$ is an isomorphism in $\A^{\nH}_{\nH},$ and using the fact that the same argument as in Remark 2.4 in \cite{G} shows that $\alpha_a$ is the counit of the adjunction $K_{\uH} \dashv R_{\uH}$, one concludes that $K_{\uH} R_{\uH}\simeq 1.$
Thus the functor $K_{\uH}$ is an equivalence of categories.
This completes the proof.
 \end{proof}

For an example,  let
 $\mathcal{V}=(\bV, \otimes, I, \sigma)$ be a braided monoidal category and
$\bH=(H, m,e,\delta, \varepsilon)$ a bialgebra in  $\mathcal{V}$.
Then
$$(H \otimes -, m \otimes -, e \otimes -, \delta \otimes -, \varepsilon \otimes -, \tau=\sigma_{H,H} \otimes -)$$
 is a bimonad on $\bV$, and it is easy to see that the category
$\bV^{\nH}_{\nH}$ of Hopf modules is just the category $\bV^{\overline{H \otimes-}}_{\underline{H \otimes -}}(\bar{\tau})=\bV^{H\otimes -}_{H \otimes -}$.

\begin{thm} {\bf Theorem.} Let $\mathcal{V}=(\bV, \otimes, I, \sigma)$ be a braided monoidal category such that idempotents split in $\bV$. Then for any bialgebra $\bH=(H, m,e,\delta, \varepsilon)$ in $\mathcal{V}$, the following are
 equivalent:
\begin{blist}
  \item $\bH$ has an antipode;
  \item the comparison functor
$$K_{H}: \bV \to \bV^{\nH}_{\nH},\quad
 V \mapsto  (H \otimes V, m \otimes V, \delta \otimes V ),\quad
f \mapsto H \otimes f,$$ is an equivalence of categories.
\end{blist}
\end{thm}

\section{Adjoints of bimonads}

This section deals with the transfer of properties of monads and comonads
to adjoint (endo-)functors. The relevance of this interplay was already
observed by Eilenberg and Moore in \cite{EM}. An effective formalism
to handle this was developed for adjunctions in 2-categories and is nicely
presented in Kelly and Street \cite{KeSt}. For our purpose we only need
this for the 2-category of categories and for convenience we recall the basic
facts of this situation here.

\begin{thm}\label{adj.nat}{\bf Adjunctions.} \em
Let $L:\A\to \B$, $R:\B\to \A$ be an adjoint pair of functors with unit and
counit $\eta, \ve$, and
$L':\A'\to \B'$, $R':\B'\to \A'$ be an adjoint pair of functors with unit and
counit $\eta'$, $\ve'$.
Given any functors $F:\A\to\A'$ and $G:\B\to \B'$,
 there is a bijection between natural transformations
$$\alpha: L'F\to GL \quad \mbox{ and }\quad \oa: FR\to R'G $$
where $\oa$ is obtained as the composite
$$FR\stackrel{\eta' {FR}}\lra R'L'FR \stackrel{R' \alpha R }\lra R'GLR
     \stackrel{R'G \ve}\lra R'G ,$$
and $\alpha$ is given as the composite
$$L'F \stackrel{L'F \eta}\lra L'FRL \stackrel{ L'\oa L}\lra L'R'GL
  \stackrel{\ve' {GF}}\lra GL.$$
In this situation, $\alpha$ and $\oa$ are called {\em mates} under the
given adjunction and this is denoted by $a \dashv \oa$. It is nicely
displayed in the diagram
 $$\xymatrix{ \A \ar[r]^L \ar[d]_F
      & \B \ar[r]^R \ar[d]^G & \A \ar[d]^F \ar@{=>}[dl]^\oa \\
  \A' \ar@{=>}[ru]^\alpha \ar[r]_{L'} & \B' \ar[r]_{R'} & \A' .} $$
Given  further

 (i) adjunctions $\tilde L: \mC \to \A$, $\tilde R :\A\to \mC$
 and $\tilde L': \mC' \to \A'$, $\tilde R' :\A'\to \mC'$
     and a functor $H:\mC\to \mC'$, or

(ii) an adjunction $L'':\A''\to \B''$, $R'':\B''\to \A''$
and functors $F':\A'\to\A''$ and $G':\B'\to \B''$,
 we get the diagram
 $$\xymatrix{
\mC\ar[r]^{\tilde L} \ar[d]_H &
 \A \ar[r]^L \ar[d]_F
      & \B \ar[r]^R \ar[d]^G & \A \ar[d]^F \ar@{=>}[dl]^\oa
       \ar[r]^{\tilde R} & \mC \ar@{=>}[dl]^\og \ar[d]^H \\
 \mC'\ar[r]_{\tilde R'} \ar@{=>}[ru]_\gamma &
 \A' \ar@{=>}[ru]_\alpha \ar[r]_{L'} \ar[d]_{F'}
   & \B' \ar[r]_{R'} \ar[d]^{G'}  & \A' \ar[d]^{F'}\ar@{=>}[dl]^\ob
    \ar[r]_{\tilde R'} & \mC' \\
 &  \A'' \ar[r]_{L''} \ar@{=>}[ru]_\beta  & \B'' \ar[r]_{R''} & \A'',} $$
yielding the mates
$$\begin{array}{lrcl}
(M1) &
L''F'F \stackrel{\beta F}\lra G'L'F \stackrel{G'\alpha}\lra G'GL & \dashv &                 F'FG  \stackrel{F'\oa}\lra F'R'G \stackrel{\ob  G} \lra R'' G'G,\\[+1mm]
(M2)& L'\tilde L'H\stackrel{L'\beta}\lra L'F\tilde L\stackrel{\alpha\tilde L}\lra
    LG \tilde L
    & \dashv & H \tilde R R \stackrel{\ob G}\lra \tilde R' R' G
            \stackrel{\tilde R'\ob}\lra \tilde R' R' G .
\end{array}$$
\end{thm}

\begin{thm}\label{prop.mate}{\bf Properties of mates.}
Let $L, L':\A\to \B$ be functors with right adjoints $R$, $R'$, respectively,
and $\alpha:L'\to L$ a natural transformation.
\begin{rlist}
\item
If $L'':\A\to \B$ is a functor  with right adjoint $R''$
and  $\beta:L''\to L'$ a  natural transformation, then
 $$\alpha \cdot \beta \, \dashv \, \ob \cdot \oa .$$
\item
If $ \tilde L:\mC\to \A$ is a functor with right adjoint $\tilde R$,
then
 $$(\alpha_ {L'}: L'\tilde L \to L \tilde L )\; \dashv \;
     (\tilde R\oa: \tilde R R \to \tilde R R').$$
\item If $ L^o:\B\to \mC$ is a functor with right adjoint $R^o$, then
$$(L^o\alpha: L^o L' \to L^o L) \; \dashv\; ( \oa {R^o}: R R^o \to R'R^o).$$
\end{rlist}
\end{thm}

\begin{proof} (i) is a special case  of \ref{adj.nat}(M1).

(ii) follows from \ref{adj.nat}(M2) by putting $\A'=\A$, $\B'=\B$,
 $\C'=\C$ and $H'=H$.

 (iii) is derived by applying \ref{adj.nat} to the diagram
$$\xymatrix{
\A\ar[r]^{L} \ar@{=}[d] &
 \B \ar[r]^{L^o} \ar@{=}[d]
      & \mC \ar[r]^{R^o} \ar@{=}[d] & \B \ar@{=}[d] \ar@{=>}[dl]^\id
       \ar[r]^{R} & \A \ar@{=>}[dl]^\oa \ar@{=}[d]  \\
 \A \ar[r]_{\tilde L'} \ar@{=>}[ru]_\alpha &
 \B \ar@{=>}[ru]_\id \ar[r]_{L^o}
   & \mC \ar[r]_{R^o}  & \B \ar[r]_{R'} & \A.} $$
\end{proof}

As observed by Eilenberg and Moore in \cite[Section 3]{EM}, for
a left adjoint endofunctor which is a monad, the right adjoint
(if it exists) is a comonad (and vice versa).
The techniques outlined above provide a convenient
and effective way to describe this transition and to prove related
properties.
Recall that for any endofunctor $L:\A\to \A$ with right adjoint $R$,
for a positive integer $n$, the powers $L^n$ have the right adjoints $R^n$.

\begin{thm}\label{adj.mon.com}{\bf Adjoints of monads and comonads.}
Let $L:\A\to \A$ be an endofunctor with right adjoint $R$.

\begin{zlist}
\item If $\uL=(L, m_L, e_L)$ is a monad, then $\oR=(R,\delta_R,\ve_R)$
  is a comonad, where $\delta_R$, $\ve_R$ are the mates of $m_L$, $e_L$
   in the diagrams
$$\xymatrix{
 \A\ar[r]^L \ar@{=}[d] & \A \ar[r]^R \ar@{=}[d]
               &\A  \ar@{=}[d] \ar@{=>}[dl]^{\ve_R}  \\
 \A \ar@{=>}[ru]_{e_L} \ar[r]_{\id} & \A \ar[r]_{\id} & \A,} \qquad
\xymatrix{
 \A\ar[r]^L \ar@{=}[d] & \A \ar[r]^R \ar@{=}[d]
               &\A  \ar@{=}[d] \ar@{=>}[dl]^{\delta_R}  \\
 \A \ar@{=>}[ru]_{m_L} \ar[r]_{HH} & \A \ar[r]_{RR} & \A.} $$

\item If $\oL=(L, \delta_L, \ve_L)$ is a comonad, then $\uR=(R,m_R,e_R)$
   is a monad where $m_R$, $e_R$ are the mates of $\delta_L$, $\ve_L$
    in the diagrams
$$\xymatrix{
 \A\ar[r]^\id \ar@{=}[d] & \A \ar[r]^\id \ar@{=}[d]
               &\A  \ar@{=}[d] \ar@{=>}[dl]^{e_R}  \\
 \A \ar@{=>}[ru]_{\ve_L} \ar[r]_{L} & \A \ar[r]_{R} & \A,} \qquad
\xymatrix{
 \A\ar[r]^{LL} \ar@{=}[d] & \A \ar[r]^{RR} \ar@{=}[d]
               &\A  \ar@{=}[d] \ar@{=>}[dl]^{m_R}  \\
 \A \ar@{=>}[ru]_{\delta_L} \ar[r]_{L} & \A \ar[r]_{R} & \A.} $$

\end{zlist}
\end{thm}
\begin{proof}
(1) Since $e_L \dashv \ve_R$ and $m_L \dashv \delta_R$, it follows from
\ref{prop.mate} (ii) and (iii) that
  $$Le_L \dashv \ve_R R, \quad e_LL\dashv R\ve_R, \quad m_LL\dashv R\delta_R,
   \quad Lm_L\dashv \delta_RR.$$
Applying \ref{prop.mate} (i) now yields
$$\begin{array}{rl}
m_L\cdot Le_L \dashv \ve_R R\cdot \delta_R, &
    m_L\cdot e_LL \dashv R\ve_R\cdot \delta_R, \\
m_L\cdot m_LL \dashv R\delta_R \cdot \delta_R, &
m_L\cdot L m_L \dashv \delta_R R\cdot \delta_R.
\end{array} $$
 Since $\uL$ is a monad we have $m_L\cdot e_LL =m_L \cdot Le_L = \id$ and
$m_L\cdot m_LL =m_L \cdot Lm_L$, implying
$$ \ve_R R\cdot \delta_R = R\ve_R \cdot \delta_R = \id \quad \mbox{ and } \quad
  R\delta_R\cdot \delta_R = \delta_R R\cdot \delta_R.$$
This shows that $\oR=(R,\delta_R,\ve_R)$ is a comonad.
\smallskip

The proof of (2) is similar.
\end{proof}

The methods under consideration also apply
to the natural transformations $LL\to LL$ which were basic for the
definition and investigation of bimonads in previous sections.
The following results were obtained in cooperation with Gabriella B\"ohm and
Tomasz Brzezi\'nski.

\begin{thm}\label{adj.dist} {\bf Adjointness and distributive laws.}
Let $L:\A\to \A$ be an endofunctor with right adjoint $R$ and
 a natural transformation $\lambda_L:LL\to LL$. This yields a mate
$\lambda_R:RR\to RR$ in the diagram
$$\xymatrix{
 \A\ar[r]^{LL} \ar@{=}[d] & \A \ar[r]^{RR} \ar@{=}[d]
               &\A  \ar@{=}[d] \ar@{=>}[dl]^{\lambda_R}  \\
 \A \ar@{=>}[ru]_{\lambda_{L}} \ar[r]_{LL} & \A \ar[r]_{RR} & \A } $$
with the following properties:

\begin{zlist}
\item $L\lambda_L \dashv \lambda_R R$ and $\lambda_LL \dashv R\lambda_R$.
\item $\lambda_L$ satisfies the Yang-Baxter equation if and only if
      $\lambda_R$ does.
\item $\lambda_L^2=\id$ if and only if $\lambda_R^2=\id$.
\item If $\uL=(L,m_L,e_L)$ is a monad and $\lambda_L$ is monad distributive,
      then $\lambda_R$ is comonad distributive for the comonad
      $\oR=(R,\delta_R,\ve_R)$.
\item If $\oL=(L,\delta_L,\ve_L)$ is a comonad and $\lambda_L$ is
      comonad distributive,
      then $\lambda_R$ is monad distributive for the comonad
      $\uR=(R,m_R,e_R)$.
\end{zlist}
 \end{thm}

 \begin{proof} (1) follows from \ref{prop.mate}, (ii) and (iii).
The remaining assertions follow by (1) and the identities in the proof
of \ref{adj.mon.com}.
\end{proof}

Recall from Definition \ref{D.1.5} that a {\em bimonad} $H$ is a monad and a comonad with compatibility conditions involving an entwining $\lambda_H:HH\to HH$.

 \begin{thm}\label{adj.bim} {\bf Adjoints of bimonads.}
Let $\bH$ be a monad $\uH=(H,m_H,e_H)$ and a comonad $\oH=(H,\delta_H,\ve_H)$
on the category $\A$. Then a right adjoint $R$ of $H$ induces a monad
$\uR=(R,m_R, e_R)$
and a comonad $\oR=(R,\delta_R, \ve_R)$ (see \ref{adj.mon.com}) and
\begin{zlist}
\item $\bH=(\uH,\oH)$ is a bimonad with entwining $\lambda_H:\uH\oH\to \oH\uH$
 if and only if $\bR=(\uR,\oR)$ is a bimonad with entwining
 $\lambda_R:\oR\uR\to \uR\oR$.

\item $\bH=(\uH,\oH)$ is a bimonad with entwining $\lambda'_H:\oH\uH\to \uH\oH$
  if and only if $\bR=(\uR,\oR)$ is a bimonad with entwining
  $\lambda'_R:\uR\oR\to \oR\uR$.

\item If $\bH=(\uH,\oH,\lambda_H)$ is a bimonad with antipode,
 then $\bR=(\uR,\oR,\lambda_R)$ is a bimonad with antipode (Hopf monad).
\end{zlist}
\end{thm}
\begin{proof}
(1) With arguments similar to those in the proof of \ref{adj.dist}
we get that $\lambda_R$ is an entwining from $\oR$ to $\uR$.
It remains to show the properties required in Definition \ref{D.1.5}.
From \ref{prop.mate}(i) we know that
$$\begin{array}{rll}
\ve_H\cdot H\ve_H \dashv e_RR\cdot e_R, & \ve_H\cdot m_H \dashv \delta_R\cdot e_R, \\
\delta_H\cdot e_H \dashv \ve_R \cdot m_R, & e_HH \cdot e_H \dashv \ve_R\cdot R\ve_R,
&
\ve_H\cdot e_H \dashv \ve_R\cdot e_R.
\end{array} $$
Thus the equalities
$$\ve_H\cdot H\ve_H = \ve_H\cdot m_H, \quad \delta_H\cdot e_H= \ve_HH\cdot e_H,
  \quad \ve_H\cdot e_H=\id$$
hold if and only if
$$e_RR\cdot e_R= \delta_R\cdot e_R, \quad \ve_R\cdot m_R=\ve_R\cdot R\ve_R, \quad \ve_R\cdot e_R=\id.$$

The transfer of the compatibility between product and coproduct \ref{D.1.18}  is seen from the corresponding diagrams
$$\xymatrix{
HH \ar[r]^-{m_H} \ar[d]_-{H \delta_H}& H \ar[r]^-{\delta_H}& HH \\
HHH \ar[rr]_-{\lambda_H H}&& HHH , \ar[u]_-{Hm_H}} \qquad
\xymatrix{
RR & \ar[l]_-{\delta_R} R & \ar[l]_-{m_R} RR
                                        \ar[d]^{\delta_RR} \\
RRR \ar[u]^{m_RR} && \ar[ll]^-{R\lambda_R} RRR . }
$$

The proof of (2) is similar.

(3) By \ref{P.1.13}, the existence of an antipode is
equivalent to the bijectivity of the morphism
$$\gamma_H = H m_H \cdot \delta_H H: HH\to HH.$$
Since $\delta_H H \dashv R m_R$ and $H m_H \dashv\delta_R R$,
  $\gamma_H$ is an isomorphism if and only if
  $\gamma_R= R m_R \cdot \delta_R R$ is an isomorphism.
\end{proof}

Functors with right (resp. left) adjoints preserve colimits (resp.
limits) and thus \ref{T.1.14} and \ref{adj.bim} imply:

\begin{thm} \label{adj.Hopf} {\bf Hopf monads with adjoints.}
 Assume the category $\A$ to admit limits or colimits.
Let $\bH=(H, m_H, e_H, \delta_H, \varepsilon_H, \lambda_H)$ be a bimonad
on $\A$ with a right adjoint bimonad
 $\bR=(R,m_R, e_R ,\delta_R,\ve_R,\lambda_R)$.
 Then the following are equivalent:
\begin{blist}
 \item  the comparison functor
   $K_{H}: \A \to \A^{\oH}_{\uH}(\lambda_H)$ is an equivalence;
\item  the comparison functor
  $K_{R}:\A \to \A^{\oR}_{\underline{R}}(\lambda_R)$
   is an equivalence;
\item $\bH$ has an antipode;
\item $\mathbf{R}$ has an antipode.
\end{blist}
\end{thm}

Finally we observe that local prebraidings are also tranferred to the
adjoint functor.

\begin{thm}\label{adj.tau} {\bf Adjointness of $\tau$-bimonads.}
Let $H$ be a monad $\uH=(H,m_H,e_H)$ and a comonad $\oH=(H,\delta_H,\ve_H)$
on the category $\A$ with a right adjoint $R$.

 If $\bH=(\uH,\oH)$ is a bimonad with double entwining $\tau_H:HH\to HH$,
 then $\bR=(\uR,\oR)$ is a bimonad with double entwining $\tau_R:RR\to RR$.

Moreover, $\tau_H$ satisfies the Yang-Baxter equation if and only if so does $\tau_R$.
\end{thm}
\begin{proof}
Most of the assertions follow immediately from \ref{adj.dist} and \ref{adj.bim}.

It remains to verify the compatibility condition \ref{D.1.2}.
For this observe that from \ref{prop.mate}(i) we get
$$\delta_H\delta_H \dashv m_Hm_H, \quad H\tau_HH \dashv R\tau_RR, \quad m_Hm_H \dashv \delta_R\delta_R, $$
and hence
$$m_Hm_H\cdot \tau_HH\cdot \delta_H\delta_H \dashv m_Rm_R\cdot R\tau_RR\cdot \delta_R\delta_R \; \mbox{ and }\; \delta_H\cdot m_M \dashv \delta_R\cdot m_R.$$
It follows that $\bH$ satifies \ref{D.1.2} if and only if so does $\bR$.
\end{proof}

\begin{thm}\label{Hopf-dual}{\bf  Dual Hopf algebras.} \em
Let $B$ be a module over a commutative ring $R$.
$B$ is a Hopf algebra if and only if the endofunctor $B\ot_R-$
on the catgeory of $R$-modules is a Hopf monad.
By \ref{adj.bim}, $B\ot_R-$ is a bimonad (with antipode) if and only if
its right adjoint functor $\Hom_R(B,-)$ is a bimonad (with antipode).
This situation is considered in more detail in \cite{BBW}.

If $B$ is finitely generated and projective as an $R$-module and
$B^*=\Hom_R(B,R)$, then  $\Hom_R(B,-)\simeq B^*\ot_R-$ and
we obtain the familiar result that $B$ is a Hopf algebra if and only if $B^*$ is.
\end{thm}

\begin{thm}\label{Hopf-group}{\bf Characterisations of groups.} \em
For any set $G$, the endofunctor
$G\times-:\mathbf{Set} \to \mathbf{Set}$ is a Hopf bimonad on the category of sets
if and only if $G$ has a group structure (e.g.\ \cite[5.20]{W-alg}).
Since the functor $\rm{Map}(G,-)$ is right adjoint to $G\times-$,
it follows from \ref{adj.Hopf} that a set $G$ is a group if and only if
the functor $\rm{Map}(G,-):\mathbf{Set} \to \mathbf{Set}$ is a Hopf monad.
\end{thm}

\bigskip

{\bf Acknowledgements.} The authors want to express their thanks
to Gabriella B\"ohm and Tomasz Brzezi\'nski for inspiring discussions
and helpful comments.
The research was started during a visit of the
first author at the Department of Mathematics at the Heinrich
Heine University of D\"usseldorf supported by the German Research
Foundation (DFG).
He is grateful to his hosts for the warm hospitality and
 to the DFG for the financial help.

\medskip

{\bf Addresses:}

{Razmadze Mathematical Institute, Tbilisi 0193, Republic of Georgia} \\
 {\small bachi@rmi.acnet.ge}

{Department of Mathematics of HHU, 40225 D\"usseldorf, Germany} \\
  {\small wisbauer@math.uni-duesseldorf.de}

\end{document}